\documentclass[12pt]{amsart}
\usepackage{amsthm}
\usepackage{amsmath}
\usepackage{amsfonts}
\usepackage{amssymb}
\usepackage[usenames]{color}
\usepackage[mathscr]{euscript}
\usepackage[hidelinks]{hyperref}
\def\natural{{\mathbb N}}
\def\zahlen{{\mathbb Z}}
\def\real{{\mathbb{R}}}
\def\dsym{\{\emptyset\}}
\def\gsym{\{\spadesuit\}}
\def\ball#1,#2.{B(#1,#2)} 
\def\clball#1,#2.{\bar B(#1,#2)} 
\def\lebmeas{{\mathcal L}^1} 
\def\eqclass#1,#2,#3.{\bigl[\bigl(#1,#2,#3\bigr)\bigr]} 
\def\lpspace#1,#2.{\setbox0=\hbox{$#2$}{\text{\normalfont L}^{\ifdim\wd0>0pt
      #2 \else Q\fi}(#1)}}
\def\lpnorm#1,#2,#3.{{\setbox0=\hbox{$#2$}\setbox1=\hbox{$#3$}\left\|#1\right\|_{\text{\normalfont L}^{\ifdim\wd0>0pt
      #2 \else Q\fi}(\ifdim\wd1>0pt #3\else \mu_p\fi)}^{\ifdim\wd0>0pt
    #2 \else Q\fi}}}
\def\llpnorm#1,#2,#3.{{\setbox0=\hbox{$#2$}\setbox1=\hbox{$#3$}\left\|#1\right\|_{\text{\normalfont L}^{\ifdim\wd0>0pt
    #2 \else Q\fi}(\ifdim\wd1>0pt #3\else \mu_p\fi)}}}
\def\indx#1.{{\text{\normalfont in}(#1)}}
\DeclareMathOperator\dom{dom}
\DeclareMathOperator\cpend{end}
\DeclareMathOperator\cpstart{str}
\DeclareMathOperator\spt{spt}
\DeclareMathOperator\xpect{E}
\DeclareMathOperator\weight{weight}
\def\hmeas#1.{\mathscr{H}^{\setbox0=\hbox{$#1\unskip$}\ifdim\wd0=0pt 1
    \else #1\fi}} 
\def\sball#1,#2,#3.{B_{#1}(#2,#3)} 
\def\scball#1,#2,#3.{\bar B_{#1}(#2,#3)} 
\DeclareMathOperator\Symbset{Symb}
\DeclareMathOperator\Wgset{Weight}
\DeclareMathOperator\ord{ord}
\DeclareMathOperator\len{len}
\DeclareMathOperator\distance{dist}
\def\bxset#1,#2,#3.{{\rm Box}\left(#1,#2,#3\right)} 
\def\cgrowth{C_{\text{\normalfont gw}}} 
\def\cgwa#1.{C_{\text{\normalfont gw},#1}} 
\def\lset#1,#2.{F(#1,#2)} 
\def\ltset#1,#2.{F_\Theta(#1 ,#2)} 
\def\munsplit{8mu} 
\def\pirange{\text{\normalfont I}_{\text{\normalfont PI}}} 
\def\neckrange{\text{\normalfont I}_{\text{\normalfont neck}}} 
\def\modulus#1,#2,#3.{\setbox1=\hbox{$#1\unskip$}\setbox4=\hbox{$#3\unskip$}\text{\normalfont
    mod}_{\ifdim\wd1=0pt P \else #1\fi}(#2; \ifdim\wd4=0pt \mu_{#2}^{(C)}\else #3\fi )}
\def\cmodulus#1,#2.{\setbox1=\hbox{$#1\unskip$}\text{\normalfont
    mod}_{\ifdim\wd1=0pt P \else #1\fi}(#2)}
\def\asycone#1,#2.{\setbox1=\hbox{$#1\unskip$}\setbox2=\hbox{$#2\unskip$}\text{\normalfont
    as-Con}(\ifdim\wd1=0pt G \else #1\fi,\ifdim\wd2=0pt \mu_G\else
  #2\fi)}
\def\wtan#1,#2.{\setbox1=\hbox{$#1\unskip$}\setbox2=\hbox{$#2\unskip$}\text{\normalfont
    w-Tan}(\ifdim\wd1=0pt X \else #1\fi,\ifdim\wd2=0pt \mu_X\else
  #2\fi)}
\newcommand{\on}{\:\mbox{\rule{0.1ex}{1.2ex}\rule{1.1ex}{0.1ex}}\:}
\numberwithin{equation}{section} 
\theoremstyle{plain}
\newtheorem{lem}[equation]{Lemma}

\newtheorem{thm}[equation]{Theorem}
\newtheorem{cor}[equation]{Corollary}
\theoremstyle{definition}
\newtheorem{defn}[equation]{Definition}
\theoremstyle{remark}

\newtheorem{rem}[equation]{Remark}
\begin{document}
\title{The Poincar\'e inequality does not improve with blow-up}
\author{Andrea Schioppa}
\address{ETHZ}
\email{andrea.schioppa@math.ethz.ch}
\thanks{The author was supported by the ``ETH Zurich Postdoctoral Fellowship Program and the Marie Curie Actions
  for People COFUND Program''}
\begin{abstract}
  For each $\beta>1$ we construct a family $F_\beta$ of metric measure
  spaces which is closed under the operation of taking weak-tangents
  (i.e.~blow-ups), and such that each element of $F_\beta$ admits a
  $(1,P)$-Poincar\'e inequality if and only if $P>\beta$.
\end{abstract}
\maketitle
\tableofcontents
\section{Introduction}
\label{sec:intro}
\subsection*{Background}
\label{subsec:background}
The abstract Poincar\'e inequality was introduced in
\cite{heinonen98} in the study of quasiconformal homeomorphisms
of metric measure spaces where points can be connected by good
families of rectifiable curves. The investigation of PI-spaces,
i.e.~metric measure spaces equipped with doubling measures and which
admit a $(1,P)$-Poincar\'e inequality for some $P\in[1,\infty)$, has
been object of intensive research.
\par One trend of investigation has focused on the infinitesimal
structure of such spaces. For example, Cheeger
\cite{cheeger99} formulated a generalization of the classical
Rademacher Differentiation Theorem which holds for PI-spaces and
showed that in such spaces the infinitesimal geometry of Lipschitz
maps is rather constrained. Moreover, this result has allowed to
formulate a notion of \emph{analytic dimension} and 
extend notions of differential geometry, like tangent and cotangent
bundles, to a large class of nonsmooth spaces which includes Carnot
groups \cite{jerison_poinc}, spaces with synthetic Ricci lower bounds
\cite{rajala_local_poinc}, some inverse limit
systems of cube complexes \cite{cheeger_inverse_poinc}, and boundaries of certain Fuchsian
buildings \cite{bourdon_pajot_poinc}. There are also more complicated examples which involve
gluing constructions \cite{hanson_poinc_nomanifold,heinonen98} and
passing to subsets \cite{mtw_carpets}. However, the
infinitesimal geometry of all these examples is rather special, in the
sense that a generic tangent/blow-up is biLipschitz equivalent to a product of
Carnot groups with an inverse limit systems of cube
complexes as in \cite{cheeger_inverse_poinc}. In general, little is
thus known about the infinitesimal structure of PI-spaces;
nevertheless, recent progress on the topic has been
achieved in \cite{cks_metric_diff}, whose results imply that a version
of \emph{metric differentiation} holds of PI-spaces, and that for a typical
blow-up $(Y,\nu)$ of a PI-space the measure $\nu$ admits a Fubini-like
representation in terms of unit speed geodesics in $Y$.
\par Another line of investigation has focused on the study of the
properties of the Poincar\'e inequality that depend on the exponent
$P$. For $\Delta>0$, a $(1,P)$-Poincar\'e inequality is stronger than
a $(1,P+\Delta)$-Poincar\'e inequality in the sense that the former
implies the latter; moreover, one can use gluing constructions to
produce examples of spaces which admit a $(1,P)$-Poincar\'e inequality
but not a $(1,P-\Delta)$-Poincar\'e inequality for some
$\Delta>0$. Intuitively, in a space admitting a $(1,P)$-Poincar\'e
inequality any pair of points can be connected by a nice family of
rectifiable curves, and the quality of these connections improves as
$P$ decreases.
\par We mention two areas of research where understanding the exponent
$P$ is important. One is the study of quasiconformal maps. For
example in \cite{kosk_mac_quasiconf_sob} it is shown that if $\varphi:X\to
Y$ is quasiconformal, where $X$ and $Y$ are metric measure spaces satisfying some
regularity assumptions (in particular $X$ is assumed to be $Q$-Ahlfors regular), if
$X$ admits a $(1,P)$-Poincar\'e inequality for $P\in[1,Q]$, so does
$Y$. However, in \cite{kosk_mac_quasiconf_sob} it is also shown that this
is not the case if $P>Q$. A second area is the study of the regularity
of minimizers and quasiminimizers of the $P$-Dirichlet energy (see for
instance \cite{kinn_sha_quasimin,kinn_mar_nonlin}); in this setting it is
usually necessary to assume a $(1,P-\Delta)$-Poincar\'e inequality for some
$\Delta>0$.
\par Given a doubling metric measure space $(X,\mu)$ we denote by
$\pirange(X,\mu)$ the largest range of exponents $P\ge 1$ such that
$(X,\mu)$ admits a $(1,P)$-Poincar\'e inequality. An open question in
analysis, even for metric spaces which can be isometrically embedded
in some Euclidean space, was whether $\pirange(X,\mu)$ is an open ray of the
form $(\beta,\infty)$. This question was answered in the affirmative in
\cite{keith_zhong_open_ended}.
\subsection*{Main Result}
\label{subsec:result}
As remarked above, as of today there is only one known class of models
for the infinitesimal geometry of PI-spaces, i.e.~biLipschitz
deformations of products of Carnot groups and inverse limit systems of
cube complexes as in \cite{cheeger_inverse_poinc}. At the same time the
first version of this preprint appeared, B.~Kleiner and the author have
found other examples~\cite{top_poinc} whose \emph{topological dimension} can
be arbitrary but whose analytic dimension is $1$. 
\par The lack of
sufficiently many examples for the infinitesimal geometry of PI-spaces
makes difficult even to formulate reasonable conjectures about the
infinitesimal geometric structure of such spaces. All the examples
mentioned above and their blow-ups at generic points
always admit a $(1,1)$-Poincar\'e inequality; while at a conference at IPAM (2013) we
learned from Le Donne of a question of Keith about whether a
$(1,P)$-Poincar\'e inequality improves to a $(1,1)$-Poincar\'e
inequality by taking tangents. Specifically, it is easy to construct
examples of $(1,p)$-PI spaces such that \emph{some} tangent does not admit a
$(1,p-\varepsilon)$-Poincar\'e inequality. For example, for $p=2$ one
can glue two copies of $\real^2$ at the origin and take on each copy
the Lebesgue measure. However, in all known examples, at \emph{a.e.~point}
\emph{all} blow-ups admit a $(1,1)$-Poincar\'e inequality.
\par In this work we answer Keith's question in the \emph{negative} and
produce new models for the infinitesimal geometry of a PI-space. In
particular, in our examples it is \emph{not} possible to improve the
Poincar\'e inequality by passing to tangents.
\begin{thm}
  \label{thm:main_result}
  There is a doubling metric space $X$ such that, for each
  $P_c\in(1,\infty)$ there exists a doubling measure $\mu_{P_c}$ on
  $X$ such that $(X,\mu_{P_c})$ and any of its weak tangents admit a
  $(1,P)$-Poincar\'e inequality if and only if $P>P_c$. The space $X$
  has Assouad-Nagata dimension $1$, and there is a Lipschitz function
  $\pi:X\to \real$ such that $(X,\mu_{P_c})$
  has a unique differentiability chart $(X,\pi)$ (i.e.~the analytic
  dimension is $1$).
\end{thm}
An interesting feature of this example is that the measures
$\{\mu_{P_c}\}_{P_c}$ can be taken \emph{mutually singular}. The existence of
$(1,1)$-Poincar\'e inequalities for mutually singular measures was
observed recently \cite{sing_poinc} in connection with the fact
that Cheeger's differentiation theorem does \emph{not} determine a canonical
measure class on a metric space. In particular, in a PI-space there
can be null sets which contain many differentiability points of a
Lipschitz function, even a common differentiability point for each
countable collection of Lipschitz functions. This is in sharp contrast
with the classical Rademacher theorem.
\par Our examples are also of interest for two different reasons. One
is that they show that there is not a strong connection between the
exponent in the Poincar\'e inequality and the underlying metric geometry of
$X$: by changing the measure class the optimal range of exponents for
which the Poincar\'e inequality holds can be arbitrarly prescribed.
\par Secondly, our examples are connected to an attempt to answer in
the \emph{negative} the question of
whether there are differentiability spaces (see~\cite{cks_metric_diff}
for details) whose infinitesimal geometry differs from that of
PI-spaces. Roughly speaking, this question asks whether a Poincar\'e
inequality is \emph{necessary}, at the infinitesimal level, to have a
Rademacher-like Theorem and a first-order calculus. The
results in~\cite{deralb, cks_metric_diff} show that differentiability spaces
share, on the infinitesimal level, similarities with PI-spaces. On the
other hand, our examples allow to move the range of exponents towards
$\infty$. The obstruction here is that degrading the range of
exponents degrades the doubling constant and so it is not possible to
get rid of the Poincar\'e inequality while keeping the measure
doubling and having first-order calculus. In a future work we
generalize the examples discussed here to overcome this problem.
\par In a forthcoming paper we also modify these examples to obtain
PI-spaces whose analytic dimension can increase by passing to
tangents. Specifically, one can have PI-spaces which are purely
$2$-unrectifiable and have analytic dimension $1$, but at generic points there
are tangents biLipschitz equivalent to $\real^2$ with the Euclidean metric.
\par Recent examples of spaces which admit $(1,P)$-Poincar\'e
inequalities but not $(1,P-\Delta)$-Poincar\'e inequalities have been
constructed in \cite{mtw_carpets,mar_spe_grad_depends}. However, such
examples are rectifiable, and so do not provide new infinitesimal
geometries. Of particular interest are the examples of
\cite{mar_spe_grad_depends} which show that the minimal $P$-weak upper
gradient depends on the choice of the exponent $P$ (i.e.~if one has a
$(1,P)$-Poincar\'e inequality but not a $(1,P-\Delta)$-Poincar\'e, the
minimal $P$-weak
upper gradient and the minimal $(P-\Delta)$-weak upper gradient
can be different). One may check that
this is not the case for our examples; this is unavoidable in
the context of taking blow-ups as discussed in~\cite{dim_blow}.
\subsection*{Overview}
\label{subsec:overview}
We observed that to produce new examples for the infinitesimal
geometry of PI-spaces one might consider an inverse limit of square
complexes where the gluing locus has $0$ 1-capacity
\cite[Example~11.13]{cheeger_inverse_poinc}. However, such examples would have
analytic and Assouad-Nagata dimension $2$, and would not give access to the
full range of exponents $P_c$. Moreover, the arguments in
\cite{cheeger_inverse_poinc} would not carry over and one would have to
resort to modulus estimates.
\par We thus decided to obtain $X$ as an asymptotic cone of a metric
graph $G$ so that the stability under blow-up would be already built
in the model. Note that one might also realize $X$ as an inverse limit
of a system of metric measure graphs, but it \emph{would not} satisfy the
same axioms as the inverse systems in \cite{cheeger_inverse_poinc}. Specifically, Axiom (2) in
\cite{cheeger_inverse_poinc}, i.e.~the requirement that simplicial projections
are open, would fail and the analysis in \cite{cheeger_inverse_poinc} would not carry over.
\par In Section~\ref{sec:exa} we construct the graph $G$ and the
corresponding measure $\mu_G$ in function of some parameters. The
choices for the weights on the measure will produce the different
measures $\mu_{P_c}$. We then make a study of the shape of balls. 
\par We point out that the definition of $G$ is somewhat technical and that
the starting point of our research were explorations of the geometry
of $G$ in C++ and Python. Specifically, it is not hard to translate
Definition~\ref{defn:G_graph} into a \texttt{Graph} class and then use Dijkstra's
shortest path algorithm to verify the results in
Subsections~\ref{subsec:walks} and~\ref{subsec:ball_box}. To help the
reader's intution we have added informal
Remarks~\ref{rem:informal_discussion1},
\ref{rem:informal_discussion2}, \ref{rem:informal_discussion_3}
and~\ref{tatiana_toro_idiot_1} to give a friendlier account of $G$.
\par In Section~\ref{sec:quasigeo} we construct \emph{good} quasigeodesics
that connect pairs of points in $G$. For convenience, we focus on the
construction of walks. To help the reader we have added an informal
discussion in Remark~\ref{rem:tatiana_toro_idiot_2}.
\par Section~\ref{sec:poinc_pf_lack} contains the technical part of
the paper. We establish modulus estimates to prove/disprove the
Poincar\'e inequality in $G$ for a given choice of $P$. In this
section we also recall the definition of modulus and a ``geometric''
characterization of the Poincar\'e inequality in terms of random
curves. 
\par Some parts of the construction of random curves are rather
technical so we provide an overview of our approach at the beginning
of Subsection~\ref{subsec:rand_curves}, and have added informal
Remarks~\ref{rem:tatiana_toro_stupid_3} and~\ref{rem:tatiana_toro_stupid_4}.
\par In Section~\ref{sec:blowup} we define asymptotic cones and
complete the proof of Theorem~\ref{thm:main_result}. In passing
information from $G$ to $X$ we take advantage of a discretization
procedure in \cite{gill_lop_disc}. 
\subsection*{Notational conventions}
\label{subsec:conventions}
We use the convention $a\approx b$ to say that $a/b,b/a\in[C^{-1},C]$
where $C$ is a universal constant; when we want to highlight $C$ we
write $a\approx_Cb$. We similarly use notations like $a\lesssim b$ and
$a\gtrsim_Cb$. In the following $C$ often denotes an unspecified universal
constant (that can change from line to line) which can be explicitly
estimated. We use the notation $E[\varphi]$ to denote the expectation
of the random variable $\varphi$. The notation $B(A,r)$ denotes a ball
of radius $R$ centred on the set $A$, i.e.~the set of points $p$ at
distance $<r$ from the set $A$.
\subsection*{Acknowledgements}
\label{subsec:thanks}
\par I want to thank E.~Le Donne for telling me about Keith's question,
and NYU/UCLA for funding the trip to IPAM.
\par I would also like to thank Bruce Kleiner for conversations that
took place at the beginning of this project, during which we discussed
variations on the constructions of~\cite{cheeger_inverse_poinc} 
in an effort to produce new examples of PI spaces.
\section{Construction of doubling graphs}
\label{sec:exa}
\subsection{Choice of parameters}
\label{subsec:par_choice}
We choose some parameters to construct the metric space $X$:
\begin{description}
\item[(P1)] An integer $N\ge2$;
\item[(P2)] A sequence of integers $\{m_k\}_k$:
  $m_k\in\left\{2,\cdots,N\right\}$ for each $k$;
\item[(P3)] Two finite sets of symbols $\Symbset_1$, $\Symbset_2$ with
  $\#\Symbset_1\ge3$ and $\#\Symbset_2\ge2$. The sets $\Symbset_1$ and
  $\Symbset_2$ share one symbol $\dsym$ which we will call the \textbf{end
    symbol}; the set $\Symbset_1$ contains another symbol $\gsym$
  which we will call the \textbf{gluing symbol}.
\end{description}
The space $X$ will be kind of self-similar in the sense that in order
to analyze its geometry we will use only a sequence of scales. We thus
introduce the scales $\sigma_k=\prod_{j=1}^k m_j$.
\par To construct the metric space $X$ we will take an asymptotic cone
of an infinite graph
$G$, see Section~\ref{sec:blowup}. For most of the paper we will work
with $G$, which is obtained as follows. We let $\Lambda$ (resp.~$\Theta$) denote the set of labels on
$\Symbset_1$ (resp.~$\Symbset_2$), i.e.~the infinite strings
$\lambda=\{\lambda(n)\}$ (resp.~$\theta=\{\theta(n)\}$) where
$\lambda(n)\in\Symbset_1$ (resp.~$\theta(n)\in\Symbset_2$) and
$\lambda(n)$ (resp.~$\theta(n)$) is eventually the end
symbol.  
\par We now regard $\real$ as a graph whose vertices are the
elements of $\zahlen$; using the scales $\sigma_k$ we associate to
each $m\in\zahlen$ an \textbf{order} $\ord(m)$ by the formula:
\begin{equation}
  \label{eq:ordef}
  \ord(m) =
  \begin{cases}
    0 &\text{if $m=0$}\\
    \max\{k: \text{$\sigma_k$ divides $|m|$}\}&\text{otherwise.}
  \end{cases}
\end{equation}
Note that if none of the $\{\sigma_k\}_k$ divides $|m|$, then by 
formula~(\ref{eq:ordef}) $\ord(m)=0$ as we convene that the $\max$ over
an empty set of natural integers is $0$.
\par We now define the graph $G$ and introduce a specific terminology
for some of its vertices.
\begin{defn}
  \label{defn:G_graph}
  Consider the graph $\real\times\Lambda\times\Theta$ and
  a vertex $v=(m,\lambda,\theta)$. Recall that we regard $\real$ as a
  graph whose vertices are the elements of $\zahlen$ and therefore
  $\real\times\Lambda\times\Theta$ is a countable union of disjoint
  graphs isomorphic to $\real$ (with vertices the elements of
  $\zahlen$ and edges of the form $[j,j+1]$ for $j\in\zahlen$). As $v$ is a vertex of
  $\real\times\Lambda\times\Theta$ recall also that $m\in\zahlen$.
  \par We say that the vertex $v$ is a \textbf{gluing
    point} of order $t$ if $\ord(m)=t>1$ and at least some symbol in
  $\{\lambda(j)\}_{j<t}$ is \textbf{not} the gluing symbol. We say that $v$
  is a \textbf{socket point} of order $t$ if $\ord(m)=t$ and
  $\lambda(j)$ is the gluing symbol for $j<t$. Note that a vertex with
  $\ord(m)=1$ is \textbf{always} a socket point.
  \par The graph $G$ is obtained from $\real\times\Lambda\times\Theta$ by
  gluing pairs of vertices
  $(m_1,\lambda_1,\theta_1),(m_2,\lambda_2,\theta_2)\in(\zahlen\times\Lambda\times\Theta)^2$
  if either on the following conditions \textbf{(Gluing)} or
  \textbf{(Socket)} holds:
  \begin{itemize}
  \item[\textbf{Gluing:}]
  \item $(m_1,\lambda_1,\theta_1),(m_2,\lambda_2,\theta_2)$ are gluing points;
  \item $m_1=m_2$ and $\theta_1=\theta_2$;
  \item $\lambda_1(j)=\lambda_2(j)$ for $j\ne\ord(m_1)$;
  \end{itemize}
  \begin{itemize}
  \item[\textbf{Socket:}]
  \item $(m_1,\lambda_1,\theta_1),(m_2,\lambda_2,\theta_2)$ are socket
    points;
  \item $m_1=m_2$;
  \item $\lambda_1(j)=\lambda_2(j)$ and $\theta_1(j)=\theta_2(j)$ for $j\ne\ord(m_1)$.    
  \end{itemize}
\end{defn}
\begin{rem}
  \label{rem:informal_discussion1}
  The previous definition of $G$ gives a precise mathematical account
  of the gluing scheme of vertices, and we used it
  to define data structures representing finite subgraphs of $G$ and
  their geodesics while we
  were exploring the connectivity properties of $G$ in C++ and
  Python. In this remark we give a more intuitive description of $G$
  to help the reader's intuition.
  \par The first step in the construction is to take countably many
  graphs isomorphic to $\real$ (where the vertices are the elements of
  $\zahlen$) and index them by pairs
  $(\lambda,\theta)\in\Lambda\times\Theta$. These graphs are just
  lines, and we can think of this union as a bunch of disjoint lines
  carrying labels and
  whose points can be represented by triples $(t,\lambda,\theta)$
  where $t$ is a ``continuous'' degree of freedom (the ``horizontal
  direction'') and $\lambda$ and
  $\theta$ are discrete degrees of freedom. In order to keep the set
  of these lines countable we impose the restriction that labels
  $\lambda$ and $\theta$ are sequences of symbols that eventually end
  in the end symbol $\dsym$.
  \par The second step is to glue the lines together to obtain a
  connected graph. Intuitively we can think of moving from a point
  $(t,\lambda,\theta)$ to a point $(s,\lambda',\theta')$, and the task
  becomes to change $t$ to $s$, $\lambda$ to $\lambda'$ and $\theta$
  to $\theta'$. Changing $t$ to $s$ does not pose a challenge as one
  can travel along the horizontal direction. To change $\lambda=\{\lambda(j)\}_{j\in\natural}$ we
  change each of the symbols $\lambda(j)$ at a time. We first focus on the case
  $j>1$; then to change $\lambda(j)$ to $\lambda'(j)$ it is 
  sufficient to reach a gluing point (or a socket point if it happens
  that the first $(j-1)$ entries of $\lambda$ are $\gsym$) traveling along the horizontal
  direction a distance $\lesssim\sigma_j$. For the
  case $j=1$ the situation is similar but we always reach a
  socket point of order $1$. Essentially the intuition is that changing
  $\lambda(j)$ is ``easy''.
  \par On the other hand, to change $\theta(j)$ to $\theta(j')$ we
  must reach a socket point $w$ of order $j$. If $j>1$ we cannot just move
  horizontally, because the $\lambda$-label of $w$ is restricted to
  have its first $(j-1)$-entries equal to the gluing symbol
  $\gsym$. Thus, socket points occur more sparsely, and unless we
  already have $\lambda(i)=\gsym$ for all $i<j$ we must first modify some
  of the labels in $\{\lambda(i)\}_{i<j}$. Essentially, the intuition
  is that changing $\theta(j)$ is ``hard'' and this will pose an
  obstruction to the existence of Poincar\'e inequalities. Note
  however, that the maximal length needed to reach a socket point of
  order $j$ is
  still $\lesssim\sigma_j$. 
  \par Finally, for $j=1$ socket points are not hard to reach as the
  restriction of their $\lambda$-label becomes vacuous. We classify
  them as ``socket points'' just because they can be used to change
  both $\lambda(1)$ and $\theta(1)$.
\end{rem}
We make $G$ a metric graph by considering the length metric where each edge has
  length $1$. Points in $G$ are then equivalence classes
  $[(t,\lambda,\theta)]$ of points $(t,\lambda,\theta)\in\real\times\Lambda\times\Theta$. The
  quotient map $\real\times\Lambda\times\Theta\to G$ will be denoted
  by $Q$. The $Q$-image of a gluing point (resp.~a socket point) will be
  called a \textbf{gluing point} (resp.~a \textbf{socket point}) of $G$. Note that the projection $\real\times\Lambda\to\real$ induces a
  $1$-Lipschitz map $\pi:G\to\real$.
\begin{rem}
  \label{rem:informal_discussion2}
  Continuing the informal discussion in
  Remark~\ref{rem:informal_discussion1}, we observe that the vertices
  of $G$ can be classified in $3$ categories. Let
  $v=[(m,\lambda,\theta)]$ be such a vertex. If $\ord(m)=j$ and if
  for $i<j$ some $\lambda(i)$ does not equal $\gsym$, then $v$
  is a gluing point of order $j$ and has valence $2\times\#\Symbset_1$. 
  If $\ord(m)=j$ and if for all $i<j$ one has $\lambda(i)=\gsym$,
  then $v$ is a socket point of order $j$ and has valence
  $2\times\#\Symbset_1\times\#\Symbset_2$. All the remaining vertices are those
  corresponding to the case $\ord(m)=0$ and have valence $2$. Finally
  note that $G$ is a graph where no edge starts and ends at the same
  point, simply because we never glue together two vertices $(m,\lambda,\theta)$
  and $(m',\lambda',\theta')$ of $\real\times\Lambda\times\Theta$ when
  $m\ne m'$. In particualr, each inclusion
  $\real\times\{\lambda\}\times\{\theta\}$ in $G$ is an isometry.
\end{rem}
To analyze the shape of balls in $G$ the following definitions
are useful.
\begin{defn}
  \label{defn:disc_log}
  To the sequence of scales $\{\sigma_k\}$ we associate the
  discretized logarithm $\lg:[0,\infty)\to\natural$ as follows:
  \begin{equation}
    \label{eq:disc_log_1}
    \lg(p) =
    \begin{cases}
      0 &\text{if $|p|<\sigma_1$}\\
      \left\{\max k: \sigma_k\le|p|\right\} &\text{otherwise.}
    \end{cases}
  \end{equation}
\end{defn}
\par Note that each vertex $v\in G$ has the form $[(k,\lambda,\theta)]$ where
$k\in\zahlen$, and $\ord(k)$ will be called the \textbf{order} of
$v$.
\begin{rem}
  \label{rem:informal_discussion_3}
  The choice of integers $\{m_k\}_k$ determines the parameters
  $\{\sigma_k\}_k$ which can be thought of as the scales at which $G$
  is analyzed. The fact that we need to consider only countably
  many scales follows from the uniform bound $2\le m_k\le N$. An
  immediate consequence of the construction is that if $v$ had order
  $k$ and $w$ has order $k'$ then $d(v,w')\ge\sigma_{\min(k,k')}$.
  \par For the
  examples considered here one can also consider the more restrictive
  case $m_k=m$ for all $k$ which gives $\sigma_k=m^k$. In this case we
  will be able to completely determine $\inf\pirange(X,\mu)$, while for
  the general case we are able only to produce upper and lower bounds
  on $\inf\pirange(X,\mu)$. In numerical experiments it seems that
  choosing the sequence $m_k$ to contain long constant subsequences
  $m_k=m_{k+1}=\ldots=m_{k+l}=i$ with different values of $i$ can
  produce tangents which have different values for the optimal
  constants in the Poincar\'e inequalities. Another effect of
  carefully choosing the sequence $\{m_k\}_k$ is an example
  $(X,\mu)$ where the analytic dimension can increase to $2$ when
  passing to tangents; this will be discussed in a forthcoming paper.
  Investigating the effects of the choice of $\{m_k\}_k$ seems an
  interesting question, in particular we leave open the question of
  whether there is a closed formula for $\inf\pirange(X,\mu)$ in terms
  of the sequence $\{m_k\}_k$ (and the other parameters in the construction).
\end{rem}
\subsection{Construction of walks}
\label{subsec:walks}
A \textbf{walk} on $G$ is a finite string on vertices and edges
$W=\{w_0\, e_1\, w_1\cdots e_l\, w_l\}$ where $w_{i-1}$ and $w_i$ are
the endpoints of $e_i$ for $1\le i\le l$. In the following we will
often suppress the edges from the notation, i.e.~simply write
$W=\{w_0\, w_1\cdots w_l\}$; we will also say that $W$ is a walk from
$w_0$ to $w_l$ and that $l$ is the length of $W$, which we will denote
by $\len W$. The starting point $\cpstart W$ of $W$ is $w_0$ and the
end point $\cpend W$ of $W$ is $w_l$. Two walks $W_1,W_2$ with $\cpend
W_1=\cpstart W_2$ can ba concatenated to obtain a walk $W_1*W_2$.
\par We say that a walk $W$ from $x$ to $y$
is \textbf{geodesic} if $\len W=d(x,y)$. This notion can be also
extended to the case in which $x$ and / or $y$ are not vertices of
$G$. In this case a \textbf{geodesic walk} from $x$ to $y$ is a
geodesic walk from a vertex $w_x$ to a vertex $w_y$ such that:
\begin{align}
  \label{eq:geowalk0}
  d(x,w_x)&<1\\
  d(y,w_y)&<1\\
  \label{eq:geowalk2}
  d(x,y)&=d(x,w_x)+\len W + d(y,w_y);
\end{align}
note that (\ref{eq:geowalk2}) implies $\len W = d(w_x,w_y)$.
A walk $W=\{w_0\,w_1\cdots w_l\}$ is \textbf{monotone increasing}
(resp.~\textbf{decreasing}) if for $0\le i\le l-1$ one has
$\pi(w_{i+1})>\pi(w_i)$ (resp.~$\pi(w_{i+1})<\pi(w_i)$).
\begin{rem}\label{tatiana_toro_idiot_1}
We have preferred to introduce walks because they are more
convenient than parametrized paths to describe the construction of
quasigeodesics and random curves that we present later in the
paper. Specifically, the following Lemmas~\ref{lem:gluing_walks},
\ref{lem:mon_lab_inc}, and~\ref{lem:mon_lab_dec} will be used to build
quasigeodesics in Section~\ref{sec:quasigeo} and to prove the
Poincar\'e inequality in Section~\ref{sec:poinc_pf_lack}.
\par In working with walks, it is important to keep track of the
labels of their vertices and edges. Recall that, except for
countably many points of $G$, the fibre $Q^{-1}(x)$ is a singleton;
the points $x$ for which $\#Q^{-1}(x)>1$ are either gluing points or
socket points. Note also that if $x$ is neither a gluing point nor a socket point, the labels
$\lambda_x\in\Lambda$ and $\theta_x\in\Theta$ are well-defined as $x=\eqclass\pi(x),\lambda,\theta.$ for 
unique $\lambda=\lambda_x$ and $\theta=\theta_x$. In particular, if $e$ is an edge, all
points in $e$, except possibly one of the vertices, have the same
labels $\lambda_e$ and $\theta_e$. 
\par On the other hand, for gluing or socket points we can still say
something about their labels. If $x$ is a
gluing point of order $k$, then $x$ is
a vertex of $G$ of the form $\eqclass\pi(x),\lambda,\theta.$ where:
$\theta$ is uniquely defined, and $\lambda(l)$ is uniquely defined for
$l\ne k$. If $x$ is a socket point of order $k$, then it is a vertex
of $G$ of the form $\eqclass\pi(x),\lambda,\theta.$ where: $\lambda(l)$ is the gluing symbol
for $l<k$, $\lambda(l)$ is uniquely defined for $l>k$, and $\theta(l)$
is uniquely defined for $l\ne k$.
Therefore, if $x$ is either a gluing point or a socket point, at most
one entry of each label $\lambda(l)$ and / or $\theta(l)$ is not uniquely
defined; in this case we will sometimes make an arbitrary choice and
still write $\lambda_x(l)$ or $\theta_x(l)$. 
\par Finally, in connection with the valence of the vertices, note that if $x$ is
a gluing point
$Q^{-1}(x)$ has cardinality $\#\Symbset_1$, and if $x$ is a socket
point $Q^{-1}(x)$ has cardinality $\#\Symbset_1\times\#\Symbset_2$.
Sometimes we will say that $\lambda$ is the $\Lambda$-label of an edge
or vertex and that $\theta$ is the $\Theta$-label of an edge or
vertex.
\end{rem}
\par In discussing walks that pass through socket points of $G$, it
will convenient to have defined a partial order on the set of labels
$\Lambda$ as one must first modify the values of the label
$\lambda$ to reach a socket pont. We say that
$\lambda<\tilde\lambda$ if there are integers $1\le k_1\le k_2$ such that:
$\lambda(j)=\tilde{\lambda}(j)$ for $j<k_1$ and $j>k_2$, and for some
$j\in[k_1,k_2]$ the entry $\tilde{\lambda}(j)$ is not the gluing symbol, and
$\lambda(j)=\gsym$ for $j\in[k_1,k_2]$.
A walk
$W=\{w_0\,e_1\,w_1\cdots e_l\,w_l\}$ is \textbf{label nondecreasing}
(resp.~\textbf{nonincreasing}) if for $1\le i \le l-1$ one has
$\lambda_{e_{i+1}}\ge \lambda_{e_i}$ (resp.~$\lambda_{e_{i+1}}\le
\lambda_{e_i}$).
\par In the following lemma we construct walks that reach a gluing (or
sometimes a socket point) moving only horizontally. They will be used
to change the value of the label $\lambda$.
\begin{lem}
  \label{lem:gluing_walks}
  Let $(p,k)\in G\times\natural$, and let $(\lambda,\theta)$ denote the labels
  of one of the edges $e$ incident to $p$. Then there is a constant
  $C$ depending only on {\normalfont\textbf {(P1)}--\textbf {(P3)}} such that there
  are monotone walks $W_+$ and $W_-$ satisfying:
  \begin{enumerate}
  \item $W_\pm$ is a walk from $p$ to $v_\pm$, where either $v_\pm$ is
    a gluing point if some $\{\lambda(j)\}_{j<k}$ is not the gluing
    symbol, or is a socket point of order $k$;
  \item $\pm(\pi(v_\pm)-\pi(p))\in[\sigma_k,C\sigma_k]$;
  \item $\len W_\pm\in[\sigma_k,C\sigma_k]$;
  \item All edges in $W_\pm$ have the same labels $(\lambda,\theta)$.
  \end{enumerate}
\end{lem}
\begin{proof}
  We just build $W_+$. Because $p$ is incident to an edge with label
  $(\lambda,\theta)$ we have $p\in
  Q(\real\times\{\lambda\}\times\{\theta\})$, and thus we can find a
  monotone increasing walk $W_0\subset
  Q(\real\times\{\lambda\}\times\{\theta\})$ which starts at $p$, has
  length $\len W_0\in[\sigma_k,2\sigma_k]$, and ends at a vertex $w_0$
  with $\ord(w_0)=0$. There is a uniform constant $C\ge 1$ such that
  the set $\real\cap[\pi(w_0),\pi(w_0)+C\sigma_k]$ contains an integer
  $t$ with $\ord(t)=k$. Let $v_+$ be the vertex of
  $Q(\real\times\{\lambda\}\times\{\theta\})$ which projects to
  $t$. Then, if all the symbols $\{\lambda(j)\}_{j\le k-1}$ equal
  $\gsym$, $v_+$ is a socket point of order $k$; otherwise $v_+$ is a
  gluing point of order $k$. Let $W_1\subset
  Q(\real\times\{\lambda\}\times\{\theta\})$ be a monotone increasing
  walk starting at $w_0$ and ending at $v_+$. Then $W_+$ is obtained
  by concatenating $W_0$ and $W_1$.
\end{proof}
\par In the following lemma we describe a walk to reach a socket
point of a given order $k$. This walk has to satisfy several technical
assumptions that we need later in the paper. Some key properties are
bounds on the length (2), the fact that the $\theta$-label is constant
(3), and restrictions (7)--(8) on the time we move in a region where a
portion of the values of the $\lambda$-label is $\gsym$.
\begin{lem}
  \label{lem:mon_lab_inc}
  Let $(p,k)\in G\times\natural$ and let $(\lambda,\theta)$ be the labels of
  an edge incident to $p$. Then there is a constant $C$
  depending on {\normalfont\textbf{(P1)}--\textbf{(P3)}} such that there are label
  nonincreasing monotone walks $W_+$ and $W_-$ satisfying:
  \begin{enumerate}
  \item $W_\pm$ is a walk from $p$ to $v_\pm$, where $v_\pm$ is a
    socket point of order $k$ such that
    $\lambda(v_\pm;l)=\lambda(p;l)$ for $l>k$;
  \item $\pm(\pi(v_\pm)-\pi(p))\in[\sigma_k,C\sigma_k]$ and $\len W_\pm\in[\sigma_k,C\sigma_k]$;
  \item The $\Theta$-label equals $\theta$  along all the edges of $W_\pm$;
  \item All the edges in $W_\pm|[0,3\sigma_k/2]$ have
    the same label $(\lambda,\theta)$;
  \item There are $(\tau_i)_{1\le i\le k-1}\subset\natural\cap[0,\len W_{\pm}]$ such that
    the map $i\mapsto\tau_i$ is strictly decreasing,
    $\tau_{k-1}\in[\frac{3\sigma_k}{2},C\sigma_k]$;
  \item The point $w_{\tau_i}$ is either a gluing point or a socket
    point of order $i$;
  \item $\len W_\pm -\tau_i\in[\sigma_i,C\sigma_i]$;
  \item Let $e_l$ be an edge of $W_\pm$; if $l\in[0,\tau_{k-1}]$,
    $\lambda_{e_l}=\lambda_{p}$; if $l\in(\tau_{i+1},\tau_i]$
    $\lambda(e_l;j)=\lambda(p;j)$ for $j\le i$ or $j>k-1$ and
    $\lambda(e_l;j)=\gsym$ for $i+1\le j\le k-1$; if
    $\lambda\in(\tau_1,\len W_\pm]$ $\lambda(e_l;j)=\gsym$ for $1\le
    j\le k-1$ and $\lambda(e_l;j)=\lambda(p_0;j)$ for $j\ge k$.
  \end{enumerate}
\end{lem}
\def\larolled#1.{\lambda^{(#1)}}
\begin{proof}
  We focus on building $W_+$ which will be built as a concatenation of
  walks $\tilde W$, $W_{k-1},W_{k-2},\ldots,W_0$.
  \par Because $p$ is incident to an edge with label
  $(\lambda,\theta)$ we have $p\in
  Q(\real\times\{\lambda\}\times\{\theta\})$, and thus we can find a
  monotone increasing walk $\tilde W\subset
  Q([\pi(p),\infty)\times\{\lambda\}\times\{\theta\})$ of length
  $\len \tilde W\in[\frac{3\sigma_k}{2},2\sigma_k]$ which starts at
  $p$ and ends at a vertex $\tilde v$ with $\ord (\tilde v) = 0$.
  \par Let $I=[\pi(\tilde v),\infty)$; in $I$ we can find a sequence
  of integers:
  \begin{equation}
    \label{eq:mon_lab_inc_p1}
    t_{k-1}\le t_{k-2}\le \cdots\le t_1\le t_0
  \end{equation}
  such that $\ord(t_i)=i$ for $i\ge 1$ and $\ord(t_0)=k$, and for some
  universal constant $C$ one has $0\le t_0-t_{k-1}\le C\sigma_k$. To
  be explicit, let $t_0$ be an integer of order $k$ in $[\pi(\tilde
  v)+\sigma_k,\pi(\tilde v)+3\sigma_k]$ and let $t_i=t_0-\sigma_i$ for
  $i\ge 1$.
  \par In the following we will let:
  \begin{equation}
    \label{eq:mon_lab_inc_p2}
    \tau_i=t_i-\pi(p),
  \end{equation}
  and we will introduce the auxiliary notation $\larolled i.$ for the
  label:
  \begin{equation}
    \label{eq:mon_lab_inc_p3}
    \larolled i.(j) =
    \begin{cases}
      \lambda(j)&\text{if $j\ge k$ or $j\le i$}\\
      \gsym&\text{otherwise.}
    \end{cases}
  \end{equation}
  \par Let $v_{k-1}$ be the vertex of
  $Q([\pi(p),\infty)\times\{\lambda\}\times\{\theta\})$ with
  $\pi(v_{k-1})=t_{k-1}$; then we let
  $W_{k-1}\subset Q([\pi(p),\infty)\times\{\lambda\}\times\{\theta\})$
  be a monotone increasing walk which starts at $\tilde v$ and ends in
  $v_{k-1}$. We let $w_{\tau_{k-1}}=v_{k-1}$ and note that $v_{k-1}$
  is either a gluing or a socket point of order $k-1$.
  \par For $i\ge 1$ the walk $W_i$ is obtained from $W_{i+1}$ as
  follows. The (backward) inductive assumption is that the last edge
  of $W_{i+1}$ has label $(\larolled i+1.,\theta)$ and that the last
  vertex $v_{i+1}$ of $W_{i+1}$ is either a gluing or a socket point
  of order $i+1$. Note that then
  $v_{i+1}\in Q([\pi(p),\infty)\times\{\larolled
  i.\}\times\{\theta\})$; we now let $v_i$ denote the vertex of
  $Q([\pi(p),\infty)\times\{\larolled i.\}\times\{\theta\})$ with
  $\pi(v_i)=t_i$. Therefore, by~(\ref{eq:mon_lab_inc_p3}) $v_i$ is
  either a gluing or a socket point of order $i$. The walk $W_i\subset
  Q([\pi(p),\infty)\times\{\larolled i.\}\times\{\theta\})$ is then
  defined as a monotone increasing walk starting at $v_{i+1}$ and
  ending in $v_i$. We then let $w_{\tau_{i}}=v_i$.
  \par We complete the construction by producing $W_0$ as follows; we
  let $\larolled 0.$ be the label such that:
  \begin{equation}
    \label{eq:mon_lab_inc_p4}
    \larolled 0.(j) =
    \begin{cases}
      \lambda(j)&\text{if $j\ge k$}\\
      \gsym&\text{otherwise.}
    \end{cases}
  \end{equation}
  We then let $v_+=v_0$ be the vertex of
  $Q([\pi(p),\infty)\times\{\larolled 0.\}\times\{\theta\})$ such that
  $\pi(v_0)=t_0$. The walk $W_0$ is then a monotone increasing walk
  joining $v_1$ to $v_0$.
  \par We now explain how each property in the statement of this
  Lemma holds:
  \begin{description}
  \item[($1$)] because $v_+=v_0$ is a socket point of order $k$ as
    $\ord(t_0)=k$ and the label $\larolled 0.$ has its first $k-1$
    entries equal to $\gsym$;
  \item[($2$)] because we have $\len \tilde W\lesssim \sigma_k$, $\len
    W_i\lesssim \sigma_i$ for $i\ge 1$ and $\len W_0\lesssim \sigma_1$;
  \item[($3$)] because the walks $\tilde W$,
    $W_{k-1},W_{k-2},\ldots,W_0$ lie in $Q(\real\times\Lambda\times\{\theta\})$;
  \item[($4$)] because of how $\tilde W$ was constructed;
  \item[($5$--$7$)] because of how the $t_i$ where chosen;
  \item[($8$)] because of how the labels $\larolled i.$ were chosen.
  \end{description}
\end{proof}
The next Lemma~\ref{lem:mon_lab_dec} is proven like
Lemma~\ref{lem:mon_lab_inc}; the proof is omitted as it looks like the
specular image of the previous one. Note that this lemma is just the
reverse situation in which we start from a socket point of order $k$
and we want to move away from it modifying the first $k$-entries of
the $\lambda$-label.
\begin{lem}
  \label{lem:mon_lab_dec}
  Let $v\in G$ be a socket point of order $k_0$ and let $\lambda$ be a
  label in $\Lambda$ such that for $k\le k_0$ one has
  $\lambda(l)=\lambda_v(l)$ for $l>k$, $l\ne k_0$. Let $\theta$ be a label
  in $\Theta$ such that $\theta_v(j)=\theta(j)$ for $j\ne k_0$.
  Then there is a constant $C$
  depending on {\normalfont\textbf{(P1)}--\textbf{(P3)}} such that there are label
  non-decreasing monotone walks $W_+$ and $W_-$  satisfying:
  \begin{enumerate}
  \item $W_\pm$ is a walk from $v$ to a vertex $p_\pm$ of order $0$
    such that $\lambda_{p_\pm}=\lambda$ and $\theta_{p_\pm}=\theta$;
  \item $\pm(\pi(p_\pm)-\pi(v))\in[\sigma_k,C\sigma_k]$ and $\len
    W_\pm\in[\sigma_k,C\sigma_k]$;
  \item All edges of $W_{\pm}$ have $\Theta$-label $\theta$;
  \item All the edges in $W_\pm|[\len W_{\pm}-\sigma_k/2,\len W_{\pm}]$ have
    the same labels;
  \item There are $(\tau_i)_{1\le i\le {k-1}}\subset\natural\cap[0,\len W_{\pm}]$ such that
    the map $i\mapsto\tau_i$ is strictly increasing, and
    $\tau_{k-1}\in[0,\len W_{\pm}-\frac{\sigma_k}{2}]$;
  \item The point $w_{\tau_i}$ is either a gluing point or a socket point of
    order $i$;
  \item $\tau_i\in[\sigma_i,C\sigma_i]$;
  \item Let $e_l$ be an edge of $W_\pm$; if $l\in[0,\tau_{1}]$,
    $\lambda{(e_l;j)}=\lambda(v;j)$ for $j\ne k_0$ and $\lambda(e_l;k_0)=\lambda(k_0)$; if $l\in(\tau_{i},\tau_{i+1}]$
    $\lambda(e_l;j)=\lambda(j)$ for $j\le i$ or $j>k-1$ and
    $\lambda(e_l;j)=\gsym$ for $i< j\le k-1$; if
    $\lambda\in(\tau_{k-1},\len W_\pm]$ $\lambda_{e_l}=\lambda$.
  \end{enumerate}
\end{lem}
\subsection{Comparison of balls and boxes}
\label{subsec:ball_box}
In the following it will be useful to replace balls by boxes
because it is easier to estimate the measure of a box;
given a Borel set $I\subset\real$, $k\in\natural\cup\{0\}$ and a finite set
  $S_1\times S_2\subset\Lambda\times\Theta$, we define the \textbf{box} $\bxset
  I,S_1\times S_2,k.$ as follows:
  \begin{equation}
    \label{eq:disc_log_2}
    \left\{[(t,\lambda,\theta)]\in G: \text{$t\in I$ and
        $\exists(\tilde\lambda,\tilde\theta)\in S_1\times S_2:$ $\forall l>k$ $(\lambda(l),\theta(l))=(\tilde\lambda(l),\tilde\theta(l))$}\right\}.
  \end{equation}
The following
lemma shows that boxes and balls are uniformly comparable.
\begin{lem}
  \label{lem:ball_box}
  Let $x=[(t,\lambda,\theta)]\in G$ and $R>0$. Let $M$ be the highest order
  of an integer $m\in[t-R,t+R]$. If $M\le\lg(2R)$ let
  $S(x,R)=\{(\lambda,\theta)\}$. If $M>\lg(2R)$ let $\Omega_M$ be the
  set of those labels $(\lambda',\theta')$
  obtained from $(\lambda,\theta)$ by making $(\lambda(M),\theta(M))$ arbitrary, and
  let $S(x,R)=\Omega_M$. Then there is a universal
  constant $C$ depending only on {\normalfont \textbf{(P1)}--\textbf{(P3)}} such
  that:
  \begin{multline}
    \label{eq:ball_box_s1}
    \bxset{\left[\pi(x)-R/2,\pi(x)+R/2\right]},\{(\lambda,\theta)\},\lg(R/C).
    \subset\clball x,R.\\
    \subset\bxset{\left[\pi(x)-R,\pi(x)+R\right]},{S(x,R)},\lg(2R)..
  \end{multline}
\end{lem}
\begin{proof}
  If $C$ is sufficiently large, using Lemmas~\ref{lem:mon_lab_inc},
  \ref{lem:mon_lab_dec} we can find, for any label $(\tilde
  \lambda,\tilde \theta)$ such that:
  \begin{equation}
    \label{eq:ball_box_p1}
    (\tilde\lambda(j),\tilde\theta(j))=(\lambda,\theta)\quad(\text{for $j>\lg(R/C)$}),
  \end{equation}
  a path of length at most $R/2$ from $x$ to a point $\tilde x$ such
  that:
  \begin{align}
    \label{eq:ball_box_p2}
    \pi(\tilde x)&=\pi(x)\\
    \label{eq:ball_box_p3}
    \tilde x&\in Q(\real\times\{\tilde \lambda\}\times\{\tilde \theta\});
  \end{align}
  this implies the inclusion:
  \begin{equation}
    \label{eq:ball_box_p4}
    \bxset{\left[\pi(x)-R/2,\pi(x)+R/2\right]},\{(\lambda,\theta)\},\lg(R/C).
    \subset\clball x,R..
  \end{equation}
  \par Let $\gamma$ be a geodesic from $x$ to $p\in\clball x,R.$; note
  that $\len\pi(\gamma)=\len\gamma$ and thus
  $\pi(\gamma(t))\in[\pi(x)-R,\pi(x)+R]$ for each
  $t\in\dom\gamma$. Therefore, if
  $(\lambda(p;k),\theta(p;k))\ne(\lambda(x;k),\theta(x;k))$, then
  $\pi(\gamma)$ passes through an integer $t_k$ of order $k$. Assume
  that $k<M$ and let $t_M\in[\pi(x)-R,\pi(x)+R]$ have order $M$; as:
  \begin{equation}
    \label{eq:ball_box_p5}
    |t_k-t_M|\ge\sigma_k,
  \end{equation}
  we conclude that $k\le\lg(2R)$. Therefore the inclusion
  \begin{equation}
    \label{eq:ball_box_p6}
    \clball x,R.
    \subset\bxset{\left[\pi(x)-R,\pi(x)+R\right]},{S(x,R)},\lg(2R).    
  \end{equation}
  follows.
\end{proof}
\subsection{Construction of measures}
\label{subsec:measures_cons}
We now turn to the construction of the measure $\mu$ on $G$. One
possibility is to take the pushforward under the quotient map
$Q:\real\times\Lambda\times\Theta\to G$ of the measure which coincides with Lebesgue
measure on each $\real\times\{\lambda\}$. For extra flexibility, in
particular to produce mutually singular measures with different values
of $\inf\pirange(X,\mu)$, we need to
choose a finite set of weights $\Wgset=\{w_s\}_{s\in\Symbset_1\cup\Symbset_2}$ subject
to the restrictions $w_s>0$ and $w_{\dsym}=1$. The restriction $w_s>0$
is needed to ensure the doubling condition, while $w_{\dsym}=1$ is
needed as our labels end eventually in $\dsym$.
\par For each
$\lambda\in\Lambda$ and $\theta\in\Theta$ we denote by $w(\lambda)$, $w(\theta)$ the associated weights:
\begin{align}
  \label{eq:label_weight1}
  w(\lambda)&=\prod_{n=1}^\infty w_{\lambda(n)},\\
  \label{eq:label_weight2}
  w(\theta)&=\prod_{n=1}^\infty w_{\theta(n)},
\end{align} where the products in (\ref{eq:label_weight1}--\ref{eq:label_weight2}) are
actually finite. We also use the notation $w((\lambda,\theta))$ for
the product $w(\lambda)w(\theta)$.
\begin{defn}
  \label{defn:measure}
  We denote by $\mu$ the measure on $G$ which is the pushforward of
  the measure on $\real\times\Lambda\times\Theta$ which coincides with
  $w((\lambda,\theta))\lebmeas$ on each
  $\real\times\{(\lambda,\theta)\}$. Note that different choices of
  the weights in $\Wgset$ will produce mutually singular measures on
  the asymptotic cone $X$, compare \cite{sing_poinc}.
\end{defn}
The next lemma provides estimates on the measures of balls and boxes.
\begin{lem}
  \label{lem:bx_measure}
  Let $S$ be a set of pairs of labels and $k\ge 1$; assume that whenever
  $(\lambda,\theta),(\lambda',\theta')\in S$ and
  $(\lambda,\theta)\ne(\lambda',\theta')$, then $(\lambda',\theta')$
  cannot be obtained from $(\lambda,\theta)$ by modifying some of the first
  $k$-entries of $\lambda$ and/or $\theta$.
  For $i=1,2$ let $\cgwa i.=\sum_{s\in\Symbset_i}w_s$; then the measure of a box is
  given by:
  \begin{equation}
    \label{eq:bx_measure_s1}
    \mu\left(\bxset
      I,S,k.\right)=\lebmeas(I)\times\cgwa 1.^k\cgwa 2.^k\sum_{(\lambda,\theta)\in
      S}\prod_{n=k+1}^\infty w(\lambda(n),\theta(n)).
  \end{equation}
  In particular, if $x=[(t,\lambda,\theta)]$:
  \begin{equation}
    \label{eq:bx_measure_s2}
    \mu\left(\clball x,R.\right)\approx R(\cgwa 1.\cgwa 2.)^{\lg
      R}\sum_{\lambda\in S(x,R)}\prod_{n=\lg R+1}^\infty w(\lambda(n),\theta(n)).
  \end{equation}
  In particular, if the $m_k$ are all equal to some $m$ and if
  $R\ge1$, we have:
  \begin{equation}
    \label{eq:bx_measure_s3}
    \mu\left(\clball x,R.\right)\approx R^{1+\log_m\cgwa
      1.+\log_m\cgwa 2.}
    \sum_{\lambda\in S(x,R)}\prod_{n=\lg R+1}^\infty w(\lambda(n),\theta(n));
  \end{equation}
  moreover, if all the weights are equal to $1$ one has:
  \begin{equation}
    \label{eq:bx_measure_s4}
    \mu\left(\clball x,R.\right)\approx R^{1+\log_m\cgwa
      1.+\log_m\cgwa 2.}.
  \end{equation}
\end{lem}
\def\tkset{T^{(k)}_{\lambda,\theta}}
\begin{proof}
  For each pair of labels $(\lambda,\theta)$ let $\tkset$ be the set
  of labels that can be obtained from $(\lambda,\theta)$ by making the
  first $k$ entries of $\lambda$ and/or $\theta$ arbitrary. We then compute as
  follows:
  \begin{equation}
    \begin{split}
      \label{eq:bx_measure_p1}
      \mu\left(
        \bxset I,S,k.
      \right) &= \sum_{(\lambda,\theta)\in S} \sum_{(\tilde
        \lambda,\tilde \theta)\in\tkset} \mu\left(
        \bxset I,S,k. \cap Q(\real\times\{\tilde\lambda\}\times\{\tilde\theta\})
      \right)\\
      &=\sum_{(\lambda,\theta)\in S} \sum_{(\tilde
        \lambda,\tilde \theta)\in\tkset}
      \lebmeas(I) w\left(
        (\tilde\lambda,\tilde\theta)
      \right)\\
      &=\lebmeas(I)\sum_{(\lambda,\theta)\in S} \sum_{(\tilde
        \lambda,\tilde \theta)\in\tkset}\prod_{n=1}^kw\left(
        (\tilde\lambda(n),\tilde\theta(n))
      \right) \cdot      
      \prod_{n=k+1}^\infty w \left(
        (\lambda(n),\theta(n))
      \right)\\
      &=\lebmeas(I)\times\cgwa 1.^k\cgwa 2.^k\sum_{(\lambda,\theta)\in
        S}\prod_{n=k+1}^\infty w(\lambda(n),\theta(n)),
    \end{split}
  \end{equation}
  which gives~(\ref{eq:bx_measure_s1}).
  \par Now~(\ref{eq:bx_measure_s2}) follows
  from~(\ref{eq:bx_measure_s1}) and Lemma~\ref{lem:ball_box} by
  observing that for any $C_0$ there is a $C(C_0)$ such that:
  \begin{equation}
    \label{eq:bx_measure_p2}
    \lg(C_0R)\le\lg R+C(C_0).
  \end{equation}
  If we assume that all the $m_k$ are equal to $m$, then the
  discretized logarithm $\lg$ is just $\log_m$ up to a bounded
  additive error, and
  hence~(\ref{eq:bx_measure_s3}), (\ref{eq:bx_measure_s4}) follow.
\end{proof}
\section{Construction of good walks}
\label{sec:quasigeo}
In this section we prove the existence of \textbf{good walks}
between points in $G$. These walks correspond to quasigeodesics which are used to build the
families of curves used to prove Poincar\'e inequalities. 
\par Let $x,y\in G$; choose labels $(\lambda_x,\theta_x)$,
$(\lambda_y,\theta_y)$ such that
$x=\eqclass\pi(x),\lambda_x,\theta_x.$, $y=\eqclass\pi(y),\lambda_y,\theta_y.$ and the
cardinality of the set:
\begin{equation}
  \label{eq:labset}
  \natural(x,y)=\left\{k:(\lambda_x(k),\theta_x(k))\ne(\lambda_y(k),\theta_y(k))\right\}
\end{equation}
is minimal.
\par In the following $C$ will denote a universal constant that can
change from line to line and that can be explicitly estimated.
\begin{defn}
  \label{defn:gwalk}
  Given $x,y\in G$ with $d(x,y)>1$ a \textbf{good walk}
  $W=\left\{w_0\,e_1\,w_1\cdots e_L\,w_L\right\}$ from $x$ to $y$ is a
  walk having the following properties:
  \begin{description}
  \item[(GW1)] $\len W\le Cd(x,y)$;
  \item[(GW2)] $d(w_0,x),d(w_L,y)\in[0,1)$;
  \item[(GW3)] for $i>0$ one has $d(w_i,x)\ge i/C$.
  \end{description}
\end{defn}
\begin{rem}
  \label{rem:tatiana_toro_idiot_2}
  Intuitively condition \textbf{(GW1)} forces the path corresponding
  to $W$ to be a
  quasigeodesic. Condition
  \textbf{(GW2)} forces $W$ to start at a vertex adjacent to an edge
  containing $x$ and end at a vertex adjacent to an edge containing $y$.
  Finally \textbf{(GW3)} forces $W$ to move away from $x$ at a linear rate.
\end{rem}
In the following we will often use the following estimate.
\begin{lem}
  \label{lem:exc_estimate}
  If $\lg d(x,y)<\max\natural(x,y)$ then for each
  $k\in\natural(x,y)\setminus\left\{\max\natural(x,y)\right\}$ one has
  $\lg d(x,y)\ge k$.
\end{lem}
\begin{proof}
  Let $w_0,w_1$ be vertices of $G$ with $\ord(w_0)\ne\ord(w_1)$, then:
  \begin{equation}
    \label{eq:exc_estimate_p1}
    d(w_0,w_1)\ge\left|\pi(w_0)-\pi(w_1)\right|\ge\sigma_{\min(\ord(w_0),\ord(w_1))}.
  \end{equation}
  Take a geodesic walk $W$ from $x$ to $y$. Then there are
  $w_{j_0},w_{j_1}\in W$ such that $w_{j_0}$ is either a gluing or a socket point of
  order $\max\natural(x,y)$ and $w_{j_1}$ is either a gluing or a socket point of order
  $k$; let $\tilde W$ be a subwalk of $W$ joining
  $w_{j_0}$ and $w_{j_1}$, and observe that:
  \begin{equation}
    \label{eq:exc_estimate_p2}
    \len(W)\ge\len(\tilde W)\ge d(w_{j_0},w_{j_1})\ge\sigma_k.
  \end{equation}
\end{proof}
\def\kmax{{k_{\text{\normalfont max}}}}
\par The following Theorem is the first part of the construction of
good walks under the additional assumption $\lg
d(x,y)\ge\kmax=\max\natural(x,y)$. Condition \textbf{(GWA2)} is just
an estimate on the length of $W$. Condition \textbf{(GWA1)} is less
transparent. It establishes what happens along $W$ as we change the
values of $\theta$ and $\lambda$ to reach $y$. If
$k\in\max\natural(x,y)$ and we need only to change $\lambda(k)$, then
after a ``critical'' vertex $w_{s(k)}$ we will always move through
edges where either $\lambda(k)=\lambda_y(k)$ or $\lambda(k)=\gsym$
(this second option occurs when we need to get closer to a socket
point of order $j>k$). If
$k\in\max\natural(x,y)$ and we need to change $\theta(k)$ (and
possibly also $\lambda(k)$), then we need to pass through a
``critical'' vertex $w_{s(k)}$ which is a socket point and after
passing through it $\theta(k)$ will remain equal to $\theta_y(k)$. An
important constraint is that the map $k\mapsto s(k)$ is monotone
increasing. Note that in the case $k\not\in\natural(x,y)$ and
$k<\kmax$ we just define $w_{s(k)}$ to be a (gluing) point along the walk so
that~(\ref{eq:gw_part1_s1}) holds.
\begin{thm}
  \label{thm:gw_part1}
  If $\lg d(x,y)\ge\kmax=\max\natural(x,y)$ there is a good walk
  $W$ from $x$ to $y$ which has the following additional properties:
  \begin{description}
  \item[(GWA1)] If $k\in\natural(x,y)$  is such that
    $\theta_x(k)=\theta_y(k)$, there is a distinguished gluing or
    socket point $w_{s(k)}$ such that each edge $e$ preceding
    $w_{s(k)}$ satisfies $\lambda_e(k)=\lambda_x(k)$, and each edge
    $e$ following $w_{s(k)}$ satisfies either
    $\lambda_e(k)=\lambda_y(k)$ or $\lambda_e(k)=\gsym$. Moreover, in
    this case all edges $e$ satisfy $\theta_e(k)=\theta_x(k)$. If
    $k\in\natural(x,y)$ is such that $\theta_x(k)\ne \theta_y(k)$,
    there is a distinguished socket point $w_{s(k)}$ such that each edge $e$ preceding
    $w_{s(k)}$ satisfies $\theta_e(k)=\theta_x(k)$ and $\lambda_e(k)=\lambda_x(k)$, and each edge
    $e$ following $w_{s(k)}$ satisfies $\theta_e(k)=\theta_y(k)$ and either
    $\lambda_e(k)=\lambda_y(k)$ or $\lambda_e(k)=\gsym$.  Moreover, the
    map $k\mapsto s(k)$ is monotone increasing and the subwalk $W_k$
    from $w_{s(k)}$ to $w_{s(k+1)}$ satisfies:
    \begin{equation}
      \label{eq:gw_part1_s1}
      \len W_k\approx\sigma_{k+1}\approx d(w_{s(k)},w_{s(k+1)});
    \end{equation}
  \item[(GWA2)] The walk $W$ satisfies:
    \begin{equation}
      \label{eq:gw_part1_s2}
      \len W\approx\max\left\{\left|\pi(x)-\pi(y)\right|,\sigma_{\kmax}\right\}.
    \end{equation}
  \end{description}
\end{thm}
\begin{proof}
  Without loss of generality we can assume $\pi(x)\le\pi(y)$. If
  $\natural(x,y)=\emptyset$ then $x,y$ lie in some
  $Q(\real\times\{\lambda\}\times\{\theta\})$ and the construction of the walk is 
  immediate. Let $w_0$ be the vertex of $G$ satisfying
  $\pi(w_0)\in[\pi(x),\pi(x)+1)$,
  $(\lambda_{w_0},\theta_{w_0})=(\lambda_x,\theta_x)$ (if the labels for $w_0$ or $x$
  are not unique, one can choose them so that equality
  holds. Note that for a non-unique label $(\lambda_p,\theta_p)$ only one entry
  $(\lambda_p(m),\theta_p(m))$ is not uniquely determined). Order the elements of
  $\natural(x,y)$ increasingly:
  \begin{equation}
    \label{eq:gw_part1_p1}
    k_0<k_1<\cdots<k_q.
  \end{equation}
  \par Now either $\theta_x(k_0)=\theta_y(k_0)$ or $\theta_x(k_0)\ne
  \theta_y(k_0)$. The goal is to
  construct a walk $W_{k_0}$ of length comparable to $\sigma_{k_0}$
  which allows to change the $k_0$-th entries of the labels. We
  build $W_{k_0}$ in two parts $W^{(-)}_{k_0}$ and $W^{(+)}_{k_0}$.
  \par We now consider the first case $\theta_x(k_0)=\theta_y(k_0)$
  which implies $\lambda_x(k_0)\ne \lambda_y(k_0)$; by
  Lemma~\ref{lem:gluing_walks} we can find a monotone increasing walk
  $W^{(-)}_{k_0}$ from $w_0$ to a gluing or a socket point
  $v^{(-)}_{k_0}$ of order $k_0$ such that:
  \begin{enumerate}
  \item $\pi(v^{(-)}_{k_0})\in\left[\pi(w_0)+\sigma_{k_0},\pi(w_0)+C\sigma_{k_0}\right]$;
  \item all edges of $W^{(-)}_{k_0}$ have the same labels $(\lambda_{w_0},\theta_{w_0})$;
  \item $\len W^{(-)}_{k_0}\in[\sigma_{k_0},C\sigma_{k_0}]$.
  \end{enumerate}
  $W^{(-)}_{k_0}$ is the first part of the walk $W_{k_0}$ and we let
  $w_{s(k_0)}=v^{(-)}_{k_0}$.
  Let $\tilde\lambda_{w_0}$ be the label which agrees with
  $\lambda_{w_0}$ except at the $k_0$-th entry
  $\tilde\lambda_{w_0}(k_0)=\lambda_y(k_0)$. The second part of the
  walk $W^{(+)}_{k_0}$ is a monotone walk of length $\len
  W^{(+)}_{k_0}\in[1,\sigma_{k_0}]$ which terminates at a vertex of
  order $0$ and whose edges have the same label $(\tilde\lambda_{w_0},\theta_{w_0})$.
  We now consider the second case $\theta_x(k_0)\ne\theta_y(k_0)$
  which is slightly more complicated. By Lemma~\ref{lem:mon_lab_inc} we can find a label-nonincreasing
  monotone walk $W^{(-)}_{k_0}$ from $w_0$ to a socket point
  $v^{(-)}_{k_0}$ such that:
  \begin{enumerate}
  \item $\pi(v^{(-)}_{k_0})\in\left[\pi(w_0)+\sigma_{k_0},\pi(w_0)+C\sigma_{k_0}\right]$.
  \item $v^{(-)}_{k_0}$ has order $k_0$ and for $l>k_0$ one has
    $(\lambda(v^{(-)}_{k_0};l),\theta(v^{(-)}_{k_0};l))=(\lambda(w_0;l),\theta(w_0;l))$.
  \item $\len W^{(-)}_{k_0}\in[\sigma_{k_0},C\sigma_{k_0}]$.
  \end{enumerate}
  $W^{(-)}_{k_0}$ is the first part of the walk $W_{k_0}$ and we let $w_{s(k_0)}=v^{(-)}_{k_0}$.
  \par By Lemma~\ref{lem:mon_lab_dec} we find a label-nondecreasing
  monotone walk $W^{(+)}_{k_0}$ from $v^{(-)}_{k_0}$ to a vertex
  $v^{(+)}_{k_0}$ of order zero satisfying:
  \begin{enumerate}
  \item $\pi(v^{(+)}_{k_0})\in\left[\pi(w_0)+\sigma_{k_0},\pi(w_0)+C\sigma_{k_0}\right]$.
  \item For $l\le k_0$ one has
    $(\lambda(v^{(+)}_{k_0};l),\theta(v^{(+)}_{k_0};l))=(\lambda(y;l),\theta(y;l))$
    and for $l>k_0$
    $(\lambda(v^{(+)}_{k_0};l),\theta(v^{(+)}_{k_0};l))=(\lambda(x;l),\theta(x;l))$. 
  \item $\len W^{(+)}_{k_0}\in[\sigma_{k_0},C\sigma_{k_0}]$.
  \item All edges of $W^{(+)}_{k_0}$ satisfy $(\lambda_e(k_0),\theta_e(k_0))=(\lambda_y(k_0),\theta_y(k_0))$.
  \end{enumerate}
  \par The construction continues by induction on $k_j$,
  i.e.~suppose we have constructed the subwalks
  $\left\{W_{k_0},\cdots,
    W_{k_j}\right\}$ which form the first part
  of $W$. The first part $W^{(-)}_{k_{j+1}}$ of $W_{k_{j+1}}$ is a label-nonincreasing monotone walk
  $W^{(-)}_{k_{j+1}}$ from $v^{(+)}_{k_j}$ to a socket point
  $v^{(-)}_{k_{j+1}}$ of order $k_{j+1}$ such that:
  \begin{enumerate}
  \item $\pi(v^{(-)}_{k_{j+1}})\in\left[\pi(v^{(+)}_{k_j})+\sigma_{k_{j+1}},\pi(v^{(+)}_{k_j})+C\sigma_{k_{j+1}}\right]$.
  \item $v^{(-)}_{k_{j+1}}$ has order $k_{j+1}$ and for $l>k_{j+1}$ one has
    $(\lambda(v^{(-)}_{k_{j+1}};l),\theta(v^{(-)}_{k_{j+1}};l))=(\lambda(v^{(+)}_{k_j};l),\theta(v^{(+)}_{k_j};l))$.
  \item $\len W^{(-)}_{k_{j+1}}\in[\sigma_{k_{j+1}},C\sigma_{k_{j+1}}]$.
  \end{enumerate}
  We then let $w_{s(k_{j+1})}=v^{(-)}_{k_{j+1}}$.
  \par By Lemma~\ref{lem:mon_lab_dec} we complete $W_{k_{j+1}}$ by finding a label-nondecreasing
  monotone walk $W^{(+)}_{k_{j+1}}$ from $v^{(-)}_{k_{j+1}}$ to a
  vertex $v^{(+)}_{k_{j+1}}$ such that:
  \begin{enumerate}
  \item $\pi(v^{(+)}_{k_{j+1}})\in\left[\pi(v^{(-)}_{k_{j+1}})+\sigma_{k_{j+1}},\pi(v^{(-)}_{k_{j+1}})+C\sigma_{k_{j+1}}\right]$.
  \item For $l\le k_{j+1}$ one has
    $(\lambda(v^{(+)}_{k_{j+1}};l),\theta(v^{(+)}_{k_{j+1}};l))=(\lambda(y;l),\theta(y;l))$
    and for $l>k_{j+1}$
    $(\lambda(v^{(+)}_{k_{j+1}};l),\theta(v^{(+)}_{k_{j+1}};l))=(\lambda(x;l),\theta(x;l))$.
  \item $\len W^{(+)}_{k_{j+1}}\in[\sigma_{k_{j+1}},C\sigma_{k_{j+1}}]$.
  \item All edges of $W^{(+)}_{k_{j+1}}$ satisfy $(\lambda_e(k_{j+1}),\theta_e(k_{j+1}))=(\lambda_y(k_{j+1}),\theta_y(k_{j+1}))$.
  \end{enumerate}
  \par When we reach $j=q$ we have constructed the first part
  $W^{(1)}$ of the walk $W$. Property \textbf{(GW3)} is satisfied
  because $W^{(1)}$ is monotone increasing and the part of
  \textbf{(GW2)} concerning $w_0$ is also satisfied; the additional
  condition \textbf{(GWA1)} is also satisfied on $W^{(1)}$, and
  needs only to be checked there because of the way in which we
  construct the second part $W^{(2)}$ of the walk.
  \par There are two cases to consider to complete the proof.
  \vskip10pt
  \noindent(\texttt{Case 1}): $\pi(v^{(+)}_{k_q})\le\pi(y)$; then
  $v^{(+)}_{k_q}$ and $y$ belong to
  $Q(\real\times\{\lambda_y\}\times\{\theta_y\})$. Therefore, $W^{(2)}$ is constructed
  by taking a geodesic walk in $Q(\real\times\{\lambda_y\}\times\{\theta_y\})$ from
  $v^{(+)}_{k_q}$ to $y$. We need only to prove \textbf{(GW1)} which
  is a consequence of \textbf{(GWA2)}:
  \begin{equation}
    \label{eq:gw_part1_p2}
    \begin{split}
      \len W &=\sum_{j=0}^q(\len W^{(-)}_{k_{j}}+\len
      W^{(+)}_{k_{j}})+\pi(y)-\pi(v^{(+)}_{k_q})\\
      &\le C\sum_{j=0}^q\sigma_{k_{j}}+\pi(y)-\pi(v^{(+)}_{k_q})\\
      &\le C\sigma_{k_q}+\pi(y)-\pi(v^{(+)}_{k_q});
    \end{split}
  \end{equation}
  however, $\pi(x)\le \pi(v^{(-)}_{k_q})\le \pi(y)$ and so
  $\sigma_{k_q}\le\pi(y)-\pi(x)$ which implies:
  \begin{equation}
    \label{eq:gw_part1_p3}
    \len W\le C(\pi(y)-\pi(x))\le Cd(x,y).
  \end{equation}
  As $\pi$ is $1$-Lipschitz and as $W$ is monotone increasing, we have $\len W\ge \pi(y)-\pi(x)$
  which completes the proof of \textbf{(GWA2)}.
  \vskip10pt
  \noindent(\texttt{Case 2}): $\pi(y)<\pi(v^{(+)}_{k_q})$; then
  $v^{(+)}_{k_q}$ and $y$ belong to $Q(\real\times\{\lambda_y\}\times\{\theta_y\})$
  and $W^{(2)}$ is constructed by taking a geodesic walk in
  $Q(\real\times\{\lambda_y\}\times\{\theta_y\})$ from $v^{(+)}_{k_q}$ to $y$; note
  that $W^{(2)}$ is monotone decreasing. Let:
  \begin{equation}
    \label{eq:gw_part1_p4}
    W^{(2)}=\left\{z_0,\cdots,z_m=v_y\right\},
  \end{equation}
  where $v_y$ is the unique vertex satisfying
  $(\lambda_y,\theta_y)=(\lambda_{v_y},\theta_{v_y})$ and $\pi(v_y)\in[\pi(y),\pi(y)+1)$. Note
  that:
  \begin{equation}
    \label{eq:gw_part1_p5}
    \pi(x)\le\pi(y)\le\pi(v_y)\le C\sigma_{k_q}+\pi(x),
  \end{equation}
  and so
  \begin{equation}
    \label{eq:gw_part1_p6}
    \begin{split}
      \sigma_{k_q}&=\len W=\sum_{j=0}^q(\len W^{(-)}_{k_{j}}+\len
      W^{(+)}_{k_{j}})+\pi(v^{(+)}_{k_q})-\pi(y)\\
      &\le C\sigma_{k_q}\le Cd(x,y),
    \end{split}
  \end{equation}
  which establishes \textbf{(GW1)}, \textbf{(GWA2)} and the part of
  \textbf{(GW2)} concerning $w_L$.
  \par If $\pi(z_m)\ge\pi(x)+\sigma_{k_{q}}/2$ then \textbf{(GW3)}
  holds for some universal constant $C$. Otherwise, let $m_0\le m$
  denote the first integer so that:
  \begin{equation}
    \label{eq:gw_part1_p7}
    \pi(z_{m_0})<\pi(x)+\sigma_{k_q}/2;
  \end{equation}
  for $\tilde m\ge m_0$ we have $d(z_{\tilde m},z_m)<\sigma_{k_q}/2$
  as $W^{(2)}$ is a monotone decreasing geodesic walk; thus:
  \begin{equation}
    \label{eq:gw_part1_p8}
    d(z_{\tilde m},x)\ge d(z_m,x)-\sigma_{k_q}/2\ge \sigma_{\kmax}/2-1,
  \end{equation}
  and so \textbf{(GW3)} holds for some universal constant $C$
  (recall that $k_q=\kmax$).
\end{proof}
In the following theorem we complete the construction of good walks by
analyzing the case $\lg d(x,y)<\kmax=\max\natural(x,y)$; essentially
this means that, as $d(x,y)$ is less than $\sigma_{\kmax}$, one is
forced to choose a particular socket or gluing point to change
$(\lambda_x(\kmax),\theta_x(\kmax))$ to
$(\lambda_y(\kmax),\theta_y(\kmax))$. Specifically, one should think about the
situation where $d(x,y)$ is insignificant next to $\sigma_{\kmax}$,
which means that geodesics from $x$ to $y$ must pass near a given
gluing or socket point. The following  condition \textbf{(GWA3)}
essentially says that we can find a gluing or socket point $u_{\kmax}$
(which must be a socket point if $\theta_x(\kmax)\ne\theta_y(\kmax)$),
then construct good walks $W_x$ and $W_y$ from $x$ to $u_{\kmax}$ and
$u_{\kmax}$ to $y$ (respectively), which satisfy the conclusions of
Theorem~\ref{thm:gw_part1}, and finally obtain $W$ concatenating $W_x$
and $W_y$.
\begin{thm}
  \label{thm:gw_part2}
  If $\lg d(x,y)<\kmax=\max\natural(x,y)$ then there is a good walk
  $W$ from $x$ to $y$ which has the following additional property:
  \begin{description}
  \item[(GWA3)] If $\theta(x,\kmax)=\theta(y,\kmax)$ there is a
    distinguished gluing or socket point $u_{\kmax}\in W$
    of order $\kmax$ such that each edge $e$ preceding $u_{\kmax}$ satisfies
    $\lambda(e;\kmax)=\lambda(x;\kmax)$ and each edge
    following $u_{\kmax}$ satisfies
    $\lambda(e;\kmax)=\lambda(y;\kmax)$. Moreover, in this case all
    edges $e$ of $W$ satisfy $\theta(e;\kmax)=\theta(x;\kmax)$. On the
    other hand, if $\theta(x,\kmax)\ne\theta(y,\kmax)$ there is a
    distinguished socket point $u_{\kmax}$ such that each edge
    preceding $u_{\kmax}$ satisfies
    $(\lambda(e;\kmax),\theta(e;\kmax))=(\lambda(x;\kmax),\theta(e;\kmax))$
    and each edge $e$ following $u_{\kmax}$ satisfies
    $(\lambda(e;\kmax),\theta(e;\kmax))=(\lambda(y;\kmax),\theta(y;\kmax))$.
    Moreover, $W$ can be
    decomposed into consecutive walks $W_x$ and $W_y$ where $W_x$ is a
    good walk from $x$ to $u_{\kmax}$ satisfying the conclusion of
    Theorem \ref{thm:gw_part1}, and $W_y$ is a good walk from
    $u_{\kmax}$ to $y$ satisfying the conclusion of Theorem \ref{thm:gw_part1}.
  \end{description}
\end{thm}
\begin{proof}\def\socks{\mathscr{U}}
  The construction in the cases $\theta_x(\kmax)=\theta_y(\kmax)$ and
  $\theta_x(\kmax)\ne\theta_y(\kmax)$ is essentially the same, and we
  thus discuss only the latter case. The properties of the labels
  $(\lambda(e;\kmax),\theta(e;\kmax))$ follow from the construction
  and Theorem~\ref{thm:gw_part1}.
  \par Take a geodesic walk $W$ from $x$ to $y$. Note that there must be a
  socket point $\tilde u\in W$ of order $\kmax$ so that:
  \begin{equation}
    \label{eq:gw_part2_p1}
    d(x,\tilde u)+d(\tilde u,y)=d(x,y);
  \end{equation}
  moreover, let $\socks$ denote the set of socket points of order
  $\kmax$ and let $u_\kmax$ be an element of $\socks$ at minimal
  distance from $x$ so that $d(x,u_\kmax)\le d(x,\tilde u)\le
  d(x,y)$. Let $k\in\natural(x,u_\kmax)$; then if $k>\kmax$ a geodesic
  walk $W$ from $x$ to $u_\kmax$ would pass through either a gluing or a socket point of
  order $k$ and by Lemma \ref{lem:exc_estimate} we would have:
  \begin{equation}
    \label{eq:gw_part2_p2}
    d(x,u_\kmax)=\len W\ge\sigma_{\kmax}>d(x,y),
  \end{equation}
  yielding a contradiction. Hence $k\le\kmax$; note that 
  $(\lambda(u_\kmax;\kmax),\theta(u_{\kmax};\kmax))$ can take any value, and hence $k<\kmax$; we
  can then take a geodesic walk from $x$ to $u_\kmax$ which must pass
  through either a gluing or a socket point of order $k$, and we apply Lemma
  \ref{lem:exc_estimate} to conclude that:
  \begin{equation}
    \label{eq:gw_part2_p3}
    d(x,u_\kmax)=\len W\ge\sigma_k. 
  \end{equation} Thus we can apply Theorem~\ref{thm:gw_part1} to
  obtain a good walk $W_x$ from $x$ to $u_\kmax$. Note that
  (\ref{eq:gw_part2_p2}) implies that
  $(\lambda(u_\kmax;l),\theta(u_{\kmax};l))=(\lambda(x;l),\theta(x;l))$ for $l>\kmax$; in particular, as
  $\kmax=\max\natural(x,y)$, if $k\in\natural(u_\kmax,y)$ we have
  $k<\kmax$. Let $W$ be a geodesic walk from $u_\kmax$ to $y$; then it
  must pass through either a gluing or a socket point of order $k$ and
  Lemma~\ref{lem:exc_estimate} implies:
  \begin{equation}
    \label{eq:gw_part2_p4}
    d(y,u_\kmax)=\len W\ge\sigma_k;
  \end{equation}
  therefore, we can apply Theorem~\ref{thm:gw_part1} to obtain a good
  walk $W_y$ from $u_\kmax$ to $y$. For later reference, we also note
  here that:
  \begin{equation}
    \label{eq:gw_part2_p5}
    d(x,u_\kmax)+d(y,u_\kmax)\in[d(x,y),3d(x,y)].
  \end{equation}
  \par The walk $W$ is obtained by concatenating $W_x$ and $W_y$ so
  that it satisfies \textbf{(GWA3)}. Property \textbf{(GW1)} follows
  observing that:
  \begin{equation}
    \label{eq:gw_part2_p6}
    \begin{split}
      \len W = \len W_x + \len W_y \le C\left(d(x,u_\kmax)+d(u_\kmax,y)\right),
    \end{split}
  \end{equation}
  and using~(\ref{eq:gw_part2_p5}) to conclude that:
  \begin{equation}
    \label{eq:gw_part2_p7}
    \len W\le Cd(x,y).
  \end{equation}
  Property \textbf{(GW2)} holds because it holds for $W_x$ and
  $W_y$. We discuss property \textbf{(GW3)} in some cases. We will
  denote by $C_1\ge2$ the constant in \textbf{(GW3)} provided by
  Theorem~\ref{thm:gw_part1}. \def\kmaxx{{k_{\text{\normalfont
          max}}^{(x)}}}\def\kmaxy{{k_{\text{\normalfont max}}^{(y)}}}
  In the following we use the notations
  $\kmaxx=\max\natural(x,u_\kmax)$ and $\kmaxy=\max\natural(u_{\kmax},y)$.
  \vskip10pt
  \noindent(\texttt{Case 1}): $\pi(x)\le\pi(u_\kmax)\le\pi(y)$.
  \vskip10pt
  \noindent(\texttt{Case 1,1}): $W_x$ and $W_y$ are both monotone. Then
  $W$ is monotone and \textbf{(GW3)} holds.
  \vskip10pt
  \noindent(\texttt{Case 1,2}): $W_x$ is not monotone and $W_y$ is
  monotone. As in Theorem~\ref{thm:gw_part1} we decompose $W_x$ is a
  first part $W^{(m)}_x$ which is monotone, and a second part
  $\{z_0,\cdots, z_m=u_\kmax\}$. Then $\len
  W_x\approx\sigma_{\kmaxx}$ and $d(z_i,x)\ge j_x(i)/C$ there
  $j_x(i)$ is the index / position of $z_i$ in the walk $W_x$, and
  $C$ is a universal constant. Let $w\in W_y$ and $j_y(w)$ denote the
  position of $w$ in $W_y$ and $j(w)$ the position in $W$. If
  $j_y(w)<2\len W_x$, then $j(w)\le 3\len W_x$ and so:
  \begin{equation}
    \label{eq:gw_part2_p8}
    d(w,x)\ge d(w,z_m=u_\kmax)\ge\frac{\len W_x}{C_1}\ge\frac{j(w)}{3C_1}.
  \end{equation}
  If $j_y(w)>2\len W_x$, then $d(w,x)\ge d(w,u_\kmax)-d(u_\kmax,x)$;
  as $W_y$ is monotone we have $d(w,u_\kmax)\ge j_y(w)$ and so:
  \begin{equation}
    \label{eq:gw_part2_p9}
    d(w,x)\ge j_y(w)-\len W_x\ge\frac{j_y(w)}{2};
  \end{equation}
  thus
  \begin{equation}
    \label{eq:gw_part2_p10}
    j(w)=j_y(w)+\len W_x\le\frac{3}{2}j_y(w),
  \end{equation}
  and so
  \begin{equation}
    \label{eq:gw_part2_p11}
    d(w,x)>\frac{j(w)}{3}.
  \end{equation}
  \vskip10pt
  \noindent(\texttt{Case 1,3}): Suppose that $W_x$ is monotone but
  $W_y$ is not. As in Theorem~\ref{thm:gw_part1} we decompose $W_y$
  in a first part $W^{(m)}_y$ which is monotone and a second part
  $\{z_0,\cdots,z_m=v_y\}$. On $W^{(m)}_y$ we obtain \textbf{(GW3)}
  as in (\texttt{Case 1,1}).
  \par Note that:
  \begin{align}
    \label{eq:gw_part2_p12}
    \len W_x&\approx\pi(u_\kmax)-\pi(x)\\
    \len W_y&\approx\len W^{(m)}_y+m\approx\sigma_{\kmaxy}\approx d(u_\kmax,y).
    \label{eq:gw_part2_p13}
  \end{align}
  Note that for each $i$ we have
  $d(z_i,u_\kmax)\ge\sigma_{\kmaxy}/C_1$. If
  $\pi(z_i)\ge\pi(u)+\sigma_{\kmaxy}/2$ we conclude that:
  \begin{equation}
    \label{eq:gw_part2_p14}
    \begin{split}
      d(z_i,x)&\ge\pi(u_\kmax)-\pi(x)+\frac{\sigma_{\kmaxy}}{2}\gtrsim
      d(x,u_\kmax)+d(u_\kmax,y)\\
      &\gtrsim_{(*)}d(x,y)
    \end{split}
  \end{equation}
  where in $(*)$ we used~(\ref{eq:gw_part2_p5}) and where the
  constant in the lower bound can be explicitly estimated in terms of
  $C_1$.
  \par Suppose that
  $\pi(z_i)\in\left[\pi(u_\kmax),\pi(u_\kmax)+\sigma_{\kmaxy}/2\right)$. Then
  any geodesic walk from $x$ to $z_i$ must pass through some socket
  point $\tilde u\in\socks$, and we would also have
  $\kmaxy\in\natural(\tilde u,y)$ so that:
  \begin{equation}
    \label{eq:gw_part2_p15}
    \begin{split}
      d(x,z_i)&\ge d(\tilde u,x)+\sigma_{\kmaxy}\ge
      d(x,u_\kmax)+\sigma_{\kmaxy}\\
      &\gtrsim d(x,u_\kmax)+d(u_\kmax,y)\gtrsim d(x,y).
    \end{split}
  \end{equation}
  The bounds~(\ref{eq:gw_part2_p14}), (\ref{eq:gw_part2_p15}) imply
  that \textbf{(GW3)} holds on $\{z_0,\cdots,z_m\}$ with a constant
  that can be computed in terms of $C_1$.
  \vskip10pt
  \noindent(\texttt{Case 1,4}): $W_x$ and $W_y$ are both not
  monotone. The argument for (\texttt{Case 1,3}) can be adapted
  noting that $d(x,u_\kmax)\approx\sigma_{\kmaxx}$.
  \vskip10pt
  \noindent(\texttt{Case 2}): $\pi(u_{\kmax})\le \pi(x)\le \pi(y)$. After
  reaching $u_\kmax$, the walk $W$ starts to move in the direction of
  increasing values of $\pi$.
  \vskip10pt
  \noindent(\texttt{Case 2,1}): $W_y$ is monotone. There is a
  $\theta>0$ depending only on \textbf{(P2)} so that
  $\sigma_{l+\theta}\ge 3\sigma_l$ for each $l$, and there is a $C_\theta$
  depending on \textbf{(P2)} so that $\sigma_{l+\theta}\le
  C_\theta\sigma_l$ for each $l$. Let $l=\lceil\lg
  d(x,u_\kmax)\rceil$ and fix $w\in W_y$. If $j(w)\le
  \sigma_{l+\theta}$ we have that any walk from $x$ to $w$ must pass
  through a socket point of order $\kmax$ and so:
  \begin{equation}
    \label{eq:gw_part2_p16}
    \begin{split}
      d(w,x)&\ge
      d(x,u_\kmax)\gtrsim\sigma_l\gtrsim\sigma_{l+\theta}\\
      &\ge j(w).
    \end{split}
  \end{equation}
  Let $j(w)>\sigma_{l+\theta}$; then $d(w,x)\ge
  d(w,u_\kmax)-d(u_\kmax,x)$; as $W_y$ is monotone, $d(w,u_\kmax)\ge
  j_y(w)$ and so:
  \begin{equation}
    \label{eq:gw_part2_p17}
    d(w,x)\ge j_y(w)-\sigma_l\gtrsim j_y(w);
  \end{equation}
  but:
  \begin{equation}
    \label{eq:gw_part2_p18}
    \begin{split}
      j(w)&=j_y(w)+\len W_x \lesssim j_y(w)+\sigma_l\\
      &\lesssim j_y(w),
    \end{split}
  \end{equation}
  and so $d(w,x)\gtrsim j(w)$ where the constant in the lower bound
  can be estimated in terms of $C_1$, $C_\theta$ and $\theta$.
  \vskip10pt
  \noindent(\texttt{Case 2,2}): $W_y$ is not monotone. We decompose
  $W_y$ as $W^{(m)}_y\cup\{z_0,\cdots,z_m=v_y\}$ and note that we
  can use (\texttt{Case 2,1}) on $W_y$. For $\{z_0,\cdots,z_m=v_y\}$
  one can adapt the argument used in (\texttt{Case 1,3}).
  \vskip10pt
  \noindent(\texttt{Case 3}): $\pi(x)\le\pi(y)\le \pi(u_{\kmax})$. This case
  can be dealt with along the lines of (\texttt{Case 2}) except in
  the case in which $W_y$ is not monotone, where a different estimate
  is required on the terminal part $\{z_0,\cdots,z_m=v_y\}$. Any walk
  from $x$ to $z_i$ must pass through socket points of orders $\kmax$
  and $\kmaxy$ so that:
  \begin{equation}
    \label{eq:gw_part2_p19}
    d(z_i,x)\ge d(x,u_\kmax)+\sigma_{\kmaxy};
  \end{equation}
  but $W_y$ is not monotone, which implies
  $\sigma_{\kmaxy}\approx\len W_y$ which gives:
  \begin{equation}
    \label{eq:gw_part2_p20}
    d(z_i,x)\gtrsim d(x,u_\kmax)+j_y(z_i);
  \end{equation}
  but $d(x,u_\kmax)\gtrsim\len W_x$ and $j(z_i)=\len W_x+j_y(w_i)$ so
  that $d(z_i,x)\gtrsim j(z_i)$.
\end{proof}
\section{The exponents for which the Poincar\'e inequality holds}
\label{sec:poinc_pf_lack}
\subsection{Geometric characterizations of the Poincar\'e inequality}
\label{subsec:poinc_pf}
The proof of the Poincar\'e inequality will involve the construction
of families of curves joining points in $G$. Overall, we have
preferred to avoid using the language of pencils of curves
employed by
\cite{semmes_finding_curves, heinonen_analysis}, and preferred a probabilistic
language. The rationale is that our construction is naturally modelled
by Markov chains, a fact that also occurrs in the examples
\cite{cheeger_inverse_poinc}. Specifically, we will deal with measurable functions defined on a
probability space which take value in the set of (Lipschitz) curves on
a metric space $X$; such maps will be called \textbf{random curves}. To a
random curve $\Gamma$ one can associate a measurable function defined
on the same probability space and which takes values in the space of
Radon measures on $X$ by $\Gamma\mapsto\|\Gamma\|$ (the length measure); such a map will be
called a \textbf{random measure}. Finally, the maps to the end and
starting points
of $\Gamma$,
$\Gamma\mapsto\cpend\Gamma$ and $\Gamma\mapsto\cpstart\Gamma$, produce
\textbf{random points} in $X$. Here for a random point we just mean a
measurable function defined on a probability space which takes values
in the set of points of $X$; alternatively, one can think of a random
point in terms of sampling points of $X$ according to some probability
measure $P$, which is the \textbf{law} of the random point. In
particular, as a random curve $\Gamma$ can be also thought in terms of
choosing a curve according to some probability law, the extremes of $\Gamma$ will be random points.
\par Finally, the support $\spt\Gamma$ of a random
curve $\Gamma$ is the set of edges that $\Gamma$ crosses in positive
measure with positive probability:
\begin{equation}
  \label{eq:support}
  \spt\Gamma=\left\{e: P_\Gamma(\|\Gamma\|(e)>0)>0\right\}.
\end{equation}
\par To disprove the Poincar\'e inequality we will use
the notion of modulus of families of curves, which we now recall.
\begin{defn}
  \label{defn:modulus}
  Let $P\ge 1$ and $A$ be a family of locally rectifiable curves in the
  metric space $X$. We say
  that a Borel function $g:X\to[0,\infty]$ is admissible for $A$ if
  for each $\gamma\in A$ one has:
  \begin{equation}
    \label{eq:modulus_1}
    \int g\,d\|\gamma\|\ge 1.
  \end{equation}
  Having fixed a background measure $\nu$ on $X$, we define the $P$-modulus of $A$,
  $\cmodulus ,{A}.$, as the infimum of:
  \begin{equation}
    \label{eq:modulus_2}
    \int g^p\,d\nu
  \end{equation}
  where $g$ ranges over the set of functions admissible for $A$. We will
  be mainly interested in modulus when $A$ is the family $A_{p,q}$ of locally rectifiable curves
  connecting two points $p$, $q$, and when $\nu$ is of the form:
  \begin{equation}
    \label{eq:modulus_3}
    \mu_{p,q}^{(C)}= \left(
      \frac{ d(p,\cdot) }
      {
        \mu\left(B(p,d(p,\cdot))\right)
      }\chi_{B(p,Cd(p,q))}
      +
      \frac{ d(q,\cdot) }
      {
        \mu\left(B(q,d(q,\cdot))\right)
      }\chi_{B(q,Cd(p,q))}
    \right)\mu,
  \end{equation}
  where $\mu$ is a doubling measure on $X$ and $C>0$. In this case we
  will use the notation $\modulus ,{p,q},.$ for the modulus of
  $A_{p,q}$ when the background measure is $\mu_{p,q}^{(C)}$.
  \par We finally recall the definition of the \textbf{Riesz potential centred
  on $p$}:
  \begin{equation}
    \label{eq:modulus_4}
    \mu_p=\frac{d(p,\cdot)}{\mu\left(
        B(p,d(p,\cdot))
        \right)}\mu.
  \end{equation}
\end{defn}
The following Theorem summarizes a geometric characterization
of $(1,P)$-Poincar\'e
inequalities. It combines results of
Heinonen-Koskela~\cite{heinonen98},
Haj{\l}asz-Koskela~\cite{haj_kosk_sobolev_met},
Keith~\cite{keith-modulus}, and Ambrosio, Di Marino and
Savar\'e~\cite{ass_modulus_duality},
and the proof is included just
for the sake of completeness. Note that we will take
Theorem~\ref{thm:pi_modulus} as the working definition of the
Poincar\'e inequality, and so we will not need to recall the usual
definition of the Poincar\'e inequality.
\begin{thm}
  \label{thm:pi_modulus}
  Let $(X,\mu)$ be a complete doubling metric measure
  space; then $P\in\pirange(X,\mu)$ if and only if one of the following
  equivalent conditions holds:
  \begin{enumerate}
  \item There is a universal constant $C$ such that for each pair of
    points $p,q\in X$ one has:
    \begin{equation}
      \label{eq:pi_modulus_s1}
      d(p,q)^{P-1}\modulus , {p,q},.\ge C;
    \end{equation}
  \item There is a universal constant $C$ such that any pair of points
    $p,q$ can be joined by a random curve $\Gamma$ satisfying:
    \begin{equation}
      \label{eq:pi_modulus_s2}
      \lpnorm{\frac{d\xpect[\|\Gamma\|]}{d\mu_{p,q}^{(C)}}},,\mu_{p,q}^{(C)}.\le Cd(p,q).
    \end{equation}
  \end{enumerate}
\end{thm}
\def\edec{_{\normalfont\text{exit}}}
\def\ldec{_{\normalfont\text{long}}}
\def\gdec{_{\normalfont\text{good}}}
\begin{proof}
    \def\lpnorm#1,#2,#3.{{\setbox0=\hbox{$#2$}\setbox1=\hbox{$#3$}\left\|#1\right\|_{\text{\normalfont L}^{\ifdim\wd0>0pt
          #2 \else Q\fi}(\ifdim\wd1>0pt #3\else \mu_{p,q}^{(C)}\fi)}^{\ifdim\wd0>0pt
        #2 \else Q\fi}}}
  The characterization of the Poincar\'e inequality in terms
  of~(\ref{eq:pi_modulus_s1}) is due to Keith~\cite{keith-modulus}, who
  built on previous results of
  Heinonen-Koskela~\cite{heinonen_analysis,heinonen98}, and Haj{\l}asz-Koskela~\cite{haj_kosk_sobolev_met}.
  \noindent\par\texttt{Step 1: (1) implies (2).}
  \noindent\par Consider the set $A$ of locally rectifiable curves
  joining $p$ to $q$; fix $M$ large to be determined later and write
  $A=A\edec\cup A\ldec\cup A\gdec$, where:
  \begin{enumerate}
  \item $A\edec$ consists of the locally rectifiable curves in $A$
    which meet $X\setminus\clball \{{p,q}\},C{d(p,q)}.$ in positive length;
  \item $A\ldec$ are the locally rectifiable curves in $A\setminus
    A\edec$ which have length $\ge Md(x,y)$;
  \item $A\gdec$ are the rectifiable curves in $A\setminus (A\edec\cup A\ldec)$.
  \end{enumerate}
  We will now fix $\mu_{p,q}^{(C)}$ as the background measure
  with respect to which we compute moduli; using the test functions
  $g\edec=0$ on $\clball \{{p,q}\},C{d(p,q)}.$ and $g\edec=\infty$
  elsewhere, and $g\ldec=M d(p,q)$ on $\clball \{{p,q}\},C{d(p,q)}.$
  and $0$ elsewhere, we see that:
  \begin{align}
    \label{eq:pi_modulus_p1}
    \cmodulus ,A\edec.&=0\\
    \label{eq:pi_modulus_p2}
    \cmodulus ,A\ldec.&\lesssim\frac{d(p,q)}{(Md(p,q))^P};
  \end{align}
  thus for $M$ sufficiently large,
  \begin{equation}
    \label{eq:pi_modulus_p3}
    d(p,q)^{P-1}\cmodulus ,A\gdec.\ge C/2.
  \end{equation}
  Instead of computing modulus on $A\gdec$ we can compute it on the
  family of measures:
  \begin{equation}
    \label{eq:pi_modulus_p4}
    \Sigma\gdec=
    \left\{
      \hmeas._\gamma:\gamma\in A\gdec
    \right\}
  \end{equation}
  Applying the main result of~\cite{ass_modulus_duality} we get a
  probability $\pi$ on $\Sigma\gdec$ such that, denoting by
  $\nu=\int_{\Sigma\gdec}\eta\,d\pi(\eta)$, we get:
  \begin{equation}
    \label{eq:eq:pi_modulus_p5}
    \llpnorm \frac{d\nu} {d\mu_{p,q}^{(C)}}, ,\mu_{p,q}^{(C)} .=\cmodulus ,\Sigma\gdec.^{-1/P};
  \end{equation}
  using (\ref{eq:pi_modulus_p3}) we conclude that:
  \begin{equation}
    \label{eq:eq:pi_modulus_p6}
    \llpnorm \frac{d\nu} {d\mu_{p,q}^{(C)}}, , \mu_{p,q}^{(C)}.\lesssim d(p,q)^{1/Q}.
  \end{equation}
  Now, to each $\eta\in\Sigma\gdec$ we can associate a unique
  unit-speed curve $\gamma:[0,\len\gamma]\to X$ such that $\hmeas
  ._\gamma=\eta$. Thus $\pi$ becomes the law of a random curve
  $\Gamma$ with $\xpect[\|\Gamma\|]=\nu$ and
  then~(\ref{eq:pi_modulus_s2}) follows
  from~(\ref{eq:eq:pi_modulus_p6}).
  \noindent\par\texttt{Step 2: (2) implies (1)}.
  \noindent\par Take a random curve $\Gamma$
  satisfying~(\ref{eq:pi_modulus_s2}) and let $g$ be admissible for
  the curves joining $p$ to $q$. Then:
  \begin{equation}
    \label{eq:eq:pi_modulus_p7}
    \begin{split}
      1&\le\xpect\left[
        \int g\,d\|\Gamma\|\right]=\int g\,d\xpect[\|\Gamma\|]\\
      &\le\llpnorm g,P,{\mu_{p,q}^{(C)}}.\llpnorm
      {        \frac{d\xpect[\|\Gamma\|]}
        {d\mu_{p,q}^{(C)}}},,{\mu_{p,q}^{(C)}}.\\
      &\le C\llpnorm g,P,{\mu_{p,q}^{(C)}}.\cdot d(p,q)^{1/Q},
    \end{split}
  \end{equation}
  and~(\ref{eq:pi_modulus_s1}) follows minimizing in $g$.
\end{proof}
\subsection{Construction of Random curves}
\label{subsec:rand_curves}
In this subsection we construct the ingredients to build the random curves
used to verify the Poincar\'e inequality. This is the subsection where
most of the technical work takes place. As we work with walks but
need to produce random curves, we
define the Lipschitz path associated to a walk as follows.
\begin{defn}
  \label{defn:walk_to_qgeo}
  To a walk $W=\{w_0\,e_1\,w_1\cdots e_l\,w_l\}$ we can canonically
  associate a $1$-Lipschitz map $\Gamma_W:[0,\len W]\to G$ by letting
  $\Gamma_W|[l,l+1]$ be a unit speed parametrization of the edge $e_l$.
\end{defn}
\par Our construction requires $3$ building blocks, which are random
curves that satisfy some constraints. These random curves will then be
concatenated in the next subsection. As an overview we offer the
following informal discussion:
\begin{itemize}
\item Theorem~\ref{thm:compression} associates to a monotone
  walk a random curve which gets ``compressed'' through a socket
  point. This situation arises when a random curve $\Gamma$ joining
  $x$ to $y$ must pass through a given socket point $\xi$. In this
  case there will be a $t_\xi$ such that $\Gamma(t_\xi)=\xi$ and so as
  $t\to t_\xi$ the random point $\Gamma(t)$ gets closer to
  $\xi$. As there is a constraint on the labels of $\xi$ and as
  $\Gamma$ is Lipschitz, the set of possible labels of the random
  point $\Gamma(t)$ will shrink as $t$ approaches $t_\xi$. 
  Intuitively, to prove a Poincar\'e
  inequality one must show that this shrinkage is not too fast,
  otherwise one cannot satisfy~(\ref{eq:pi_modulus_s2}).
\item Theorem~\ref{thm:parallel_transport} associates to a monotone walk $W_0$
  a random curve which moves ``parallel'' to $W_0$. This
  situation arises when we have a random curve $\Gamma$ which can take
  a finite set of values which are all \emph{lifts} (compare
  Definition~\ref{defn:walk_lift}) of a given curve.
\item Theorem~\ref{thm:expansion} which explains how to ``expand'' a
  random curve so that as $t$ increases the set of possible labels for
  $\Gamma(t)$ increases. Note that this situation is already familiar
  in the classical Poincar\'e inequality. For example, consider a random
  curve $\Gamma$ joining $x$ and $y$ with $\dom\Gamma=[0,L]$ which is used to verify a Poincar\'e
  inequality by proving~(\ref{eq:pi_modulus_s2}). One expects that
  as $t\to L/2$ the random point $\Gamma(t)$ can take a broader set of
  values, leading to a more diffused probability measure. On the other hand, as
  $t\to 0$ (resp.~$t\to L$) one expects that the probability
  associated to $\Gamma(t)$ concentrates on $x$ (resp.~$y$).
\end{itemize}
We now define a notion of lift for walks used in the subsequent
constructions. The idea is that given points $w_0$ and $w'_0$
satisfying $\pi(w'_0)=\pi(w_0)$ we can canonically lift a walk
starting at $w_0$ to a walk starting at $w'_0$.
\begin{defn}
  \label{defn:walk_lift}
  Let $W=\{w_0\,e_1\,w_1\cdots e_l\,w_l\}$ and $w'_0$ a point such
  that $\pi(w'_0)=\pi(w_0)$. We construct a new walk
  $\{w'_0\,e'_1\,w'_1\cdots e'_l\,w'_l\}$ as follows. The vertex
  $w'_{i+1}$ is adjacent to $w'_i$ and is determined as follows. If
  $w'_i$ is not a socket point the requirement
  $\pi(w'_{i+1})=\pi(w_{i+1})$ uniquely determines
  $w'_{i+1}$. Otherwise, assume that $w'_i$ is a socket point of order
  $k$ and let $e'_{i+1}$ denote the edge between $w'_i$ and $w'_{i+1}$. We
  require that $\lambda(e'_{i+1};k)=\lambda(e_{i+1};k)$ and
  $\theta(e'_{i+1};k)=\theta(e_{i+1};k)$ for all $k$. We say that
  $W'$ is the \textbf{lift of $W$ starting at $w'_0$} and we will
  denote it by $w'_0\cdot W$.
\end{defn}
We now add some auxiliary definitions used in the constructions,
e.g.~when concatenating random curves. The idea is that when we need
to concatenate a random curve $\Gamma_0$ to a random curve $\Gamma_1$
we need the probability measures associated to $\cpend\Gamma_0$ and
$\cpstart\Gamma_1$ to be compatible. We thus introduce canonical
probabilities on subsets of $\pi^{-1}(s)$ (where $s\in\zahlen$) 
determined by constraints on
$\lambda$ and $\theta$. 
\begin{defn}
  \label{defn:label_set}
  Let $p\in G$ a vertex with $\ord(p)=0$ and $k\in\natural$. Let
  $\lset p,k.$ denote the set of those $p'\in G$ satisfying
  $\pi(p')=\pi(p)$ and $(\lambda_{p'}(l),\theta_{p'}(l))=(\lambda_p(l),\theta_p(l))$ for $l>k$.
  For $k=0$ we let $\lset p,0.=\{p\}$. To $\lset p,k.$ we can
  associate a canonical probability measure $P$, which can also be
  thought of as the law of a
  random point in $\lset p,k.$. The probability $P$ satisfies:
  \begin{equation}
    \label{eq:label_set_1}
    \frac{P(p')}{P(p'')}=\frac{w((\lambda_{p'},\theta_{p'}))}{w((\lambda_{p''},\theta_{p''}))}\quad(\forall
    p,p'\in\lset p,k.).
  \end{equation}
  For $p'\in\lset p,k.$ denote by $s(p')$ the finite string of pairs
  $\{(\lambda_p(j),\theta_p(j))\}_{j\le k}$; then:
  \begin{equation}
    \label{eq:label_set_2}
    P(p')=(\cgwa 1.\cgwa 2.)^{-k}w(s(p')).
  \end{equation} 
  Given $\lset p_0,k.$, $\lset p_1,k.$ we define a canonical map
  $\tau:\lset p_0,k.\to\lset p_1,k.$ so that $\tau(p'_0)$ is the
  unique point $p'_1\in\lset p_1,k.$ such that $s(p'_0)=s(p'_1)$. Note
  that $\tau_\#P_0=P_1$.
  \par Let $p\in G$ a vertex and $k\in\natural$. We denote by $\ltset
  p,k.$ the set of those $p'\in G$ satisfying $\pi(p')=\pi(p)$,
  $\lambda_{p'}=\lambda_{p}$ and $\theta_{p'}(l)=\theta_p(l)$ for
  $l>k$. As above, to $\ltset p,k.$ we associate a canonical
  probability $P$ by requiring:
  \begin{equation}
    \label{eq:label_set_3}
    \frac{P(p')}{P(p'')}=\frac{w(\theta_{p'})}{w(\theta_{p''})}\quad(\forall
    p,p'\in\ltset p,k.).
  \end{equation}
\end{defn}
\par We now present the construction of a random curve which goes
through a socket point $\xi$ in $G$ if one has a walk that passes
through $\xi$. In the following, given a walk $W=\{w_0\,e_1\,w_1\cdots
e_l\,w_l\}$ we denote by $W^{-1}$ the reversed walk
$\{w_l,e_l,w_{l-1},\cdots,e_1,w_0\}$.
\def\jcut{J_{\text{\normalfont cut}}}
\def\lset#1,#2.{\setbox2=\hbox{$#2$\unskip}F(#1;\ifdim\wd2>0pt
  #2\else k-\jcut\fi)}
\def\ltset#1,#2.{\setbox2=\hbox{$#2$\unskip}F_\Theta(#1;\ifdim\wd2>0pt
  #2\else k-\jcut\fi)}
\begin{thm}
  \label{thm:compression}
  Let $W_0$ be a monotone walk. Let $p_0=\cpstart W_0$, $\xi=\cpend
  W_0$. Assume that:
  \begin{description}
  \item[(H1)] $\ord(p_0)=0$ and $\xi$ is a socket point of
    order $K\ge k$;
  \item[(H2)] $\len W_0\in[\sigma_k,C_0\sigma_k]$ and all edges of
    $W_0$ have the same $\Theta$-label $\theta$;
  \item[(H3)] There are $(\tau_i)_{1\le i\le k-1}\subset\natural\cap[0,\len W_0]$ such that
    the map $i\mapsto\tau_i$ is strictly decreasing,
    $\len W_0-\tau_i\in[\sigma_i,C_0\sigma_i]$, $w_{\tau_i}$ is either
    a gluing or a socket point of
    order $i$, and if $l\ge\tau_i+1$ one has $\lambda(e_l;j)=\gsym$ for $i\le j\le k-1$;
  \item[(H4)] If $w_s\in W$ satisfies $\ord (w_s)\ge k$, then $\lambda_{e_s}=\lambda_{e_{s+1}}$;
  \item[(H5)] For an edge $e_t$ of $W_0$ one has the
    following: if $t\in[1,\tau_{k-1}]$ then $\lambda_{e_t}=\lambda_{p_0}$;
    if $t\in(\tau_{i+1},\tau_i]$ then $\lambda_{e_t}(l)=\lambda_{p_0}(l)$ for
    $l<i$ or $l\ge k$; if $t\in[\tau_1,\len W_0]$ then
    $\lambda_{e_t}(l)=\lambda_{p_0}(l)$ for $l\ge k$.
  \end{description}
  Fix $\jcut\in\natural\cup\{0\}$ and let $P_0$ be the canonical probability on $\lset p_0,.$. Construct a
  random curve $\Gamma$ as follows: choose $p'_0\in\lset p_0,.$
  according to $P_0$ and let $\Gamma=\Gamma_{p'_0\cdot W_0}$. Then:
  \begin{description}
  \item[(C1)] $\cpend\Gamma$ has law $P_1$, where $P_1$ is the
    canonical probability on $\ltset\xi,.$;
  \item[(C2)] $\spt\Gamma\subset \ball\Gamma_W,C_1\sigma_{k-\jcut}.$;
  \item[(C3)] To each $e\in\spt\Gamma$ there is associated a unique
    $\indx e.$ such that $\pi(e)=\pi(e_{\indx e.})$, where
    $e_{\indx e.}$ is the $\indx e.$-th edge of $W_0$, and one has:
    \begin{multline}
      \label{eq:compression_s1}
      \frac{d\xpect[\|\Gamma\|]}{d\mu}|e\approx_{C_1}\cgwa 1.^{-\lg(
        \len W_0-\indx e.)}\cgwa 2.^{-k+\jcut}\\\times w_{\gsym}^{-k+\lg(
        \len W_0-\indx e.)}
      \prod_{j=k}^\infty w(\lambda(e;j),\theta(e;j))^{-1},
    \end{multline}
    where $C_1$ depends on $\jcut$, $C_0$,
    {\normalfont \textbf{(P1)}--\textbf{(P3)}} and $\Wgset$.
  \end{description}
\end{thm}
\begin{rem}
  \label{rem:tatiana_toro_stupid_3}
  \par While the hypotheses \textbf{(H1)} and \textbf{(H2)} are clear,
  we offer more motivation for \textbf{(H3)}--\textbf{(H5)}. 
  Condition \textbf{(H3)} is an assumption on how fast the labels of
  $\lambda_e$ of the edges of $W_0$ approach $\lambda_\xi$. The point
  is that the entries of $\lambda_e$ are switched to $\gsym$ 
  in reverse order, from $k-1$ to $1$, and that switching
  $\lambda_e(i)$ occurs at a distance from $\xi$ comparable to
  $\sigma_i$. In condition \textbf{(H4)} we assume that if we pass
  through a gluing or socket point of order $>k$ we do not use it to
  change $\lambda$. Finally \textbf{(H5)} is a consistency condition
  for \textbf{(H3)}: we change as few labels as possible and after
  switching $\lambda(i)$ at $w_{\tau_i}$ we do not switch the value of
  $\lambda(i)$ again. Moreover, labels $\lambda(l)$ are never changed
  for $l\ge k$.
  \par Concerning the conclusions, we point out that \textbf{(C3)} is
  the technical estimate quantifying that the ``compression'' of
  $\Gamma$ is not too fast. This plays a crucial role in establishing
  the Poincar\'e inequality. Finally, note that $\jcut$ is an integer parameter chosen
  for convenience, i.e.~to create some ``space'' between the length of $W$
  and the maximum order of the entries of $\lambda$ and $\theta$ that
  can differ from the corresponding values in $\lambda_{p_0}$ and $\theta_{p_0}$.
\end{rem}
\begin{proof}
  We prove \textbf{(C1)}. Let $\xi'=\cpend(p'_0\cdot W_0)$; we use the
  notation $w_t,e_t$ for the vertices, respectively the edges of
  $W_0$; we use the notation $w'_t,e'_t$ for the corresponding edges
  and vertices of $p'_0\cdot W_0$. We note that if $t\ge\tau_i+1$
  \textbf{(H3)} implies that $\lambda(e'_t;l)=\gsym$ for $i\le l\le
  k-1$. We thus conclude that $\lambda(\xi';l)=\gsym$ for $l\le k-1$;
  for $l\ge k$ the label $\lambda_{e'_t}$ coincides with that of
  $\lambda_{e_t}$ and so we conclude that
  $\lambda(\xi';l)=\lambda(\xi;l)$ for $l\ge k$.
  Therefore, $\xi'$ is a socket point of order $K$. By
  \textbf{(H2)} all edges of $W_0$ have the same label $\theta$, and
  this implies that all edges of $p'_0\cdot W_0$ have the same label
  $\theta_{p'_0}$. As $\pi(\xi')=\pi(\xi)$, we conclude that $\xi'$
  is the point of $\ltset \xi,.$ with label $\theta_{p'_0}$ and thus
  \textbf{(C1)} follows.
  \par We now prove \textbf{(C2)}. Note that the $i$-th vertices
  $w_i$, $w'_i$ of $W_0$ and $p'_0\cdot W_0$ have
  $\pi(w_i)=\pi(w'_i)$, and the labels $(\lambda(w_i),\theta(w_i))$, $(\lambda(w'_i),\theta(w'_i))$ and can differ
  only in the first $k-\jcut$ entries. Hence \textbf{(C2)} follows
  from Lemma~\ref{lem:ball_box}.
  \par We now prove \textbf{(C3)}. First let $e\in\spt\Gamma$ and
  assume that $e=e'_l\in p'_0\cdot W_0$, $e=e''_{\tilde l}\in
  p''_0\cdot W_0$. As the path $W_0$ is monotone, $l=\tilde l$ and
  there is a unique edge $e_{s}$ of $W_0$ such that
  $\pi(e)=\pi(e_s)$. We can thus associate to $e$ the unique integer
  $\indx e.=s$. We now turn to the proof
  of~(\ref{eq:compression_s1}). For $p'_0\in\Lambda(p_0,k-\jcut)$ we
  will denote by $e(p'_0;l)$ the $l$-th edge of $p'_0\cdot W_0$.
  \par We
  now fix $e\in\spt\Gamma$ and assume that $\indx
  e.=s$. We first consider the case $s\in [1,\tau_{k-1}]$. Then by
  \textbf{(H5)} there is a unique $p'_0\in\lset p_0,.$ such that $e$
  is the $s$-th edge of $p'_0\cdot W_0$. In this case by \textbf{(H2)}--\textbf{(H3)}
  $\lg(\len W_0-\indx e.)$ is comparable to $k$ up to a multiplicative
  constant depending on $C_0$. Assume now that
  $s\in(\tau_i,\tau_{i+1}]$; then $e$ is the $s$-the edge of
  $p'_0\cdot W_0$ if and only if:
  \begin{align}
    \label{eq:compression_p6_1}
    \theta_{p'_0}&=\theta_{e_s}\\
    \label{eq:compression_p6_2}
    \lambda(p'_0;j)&=\lambda(e;j)\quad(1\le j<i);    
  \end{align}
  note also that in this case $\lg(\len W_0-\indx e.)$ is comparable
  to $i$. Finally by \textbf{(H3)} if $s\in[\tau_1,\len W_0]$ $e$ is the $s$-th egdge
  of $p'_0\cdot W_0$ whenever $p'_0\in\lset p_0,.$ satisfies
  $\theta_{p'_0}=\theta_e$. Note that in this case $\lg(\len W_0-\indx
  e.)$ is comparable to $1$. We can now put all this information together:
  \begin{equation}
    \label{eq:compression_p8}
    \begin{split}
      \weight\left(\xpect\left[\|\Gamma\|\right];e\right)&=\sum\left\{
        P_0(p'_0): e(p'_0;\indx e.) = e, p'_0\in\lset p'_0,k-\jcut.
      \right\}
      \\
      &\approx\sum\biggl\{P_0(p'_0):
      \theta_{p'_0}=\theta_e, \lambda\left( e(p'_0;\indx e.); j
      \right) = \lambda(e;j)\\
      &\mskip 125mu\text{ for $j<\lg (\len W-\indx e.)$}
      \biggr\}\\
      &\approx \cgwa
      2.^{-k+\jcut}\prod_{j=1}^{k-\jcut}w\left(\theta_e(j)\right)
      \cgwa 1.^{-\lg(\len W_0-\indx e.)}\\&\mskip \munsplit \times\prod_{j=1}^{\lg(\len
        W_0-\indx e.)}w(\lambda(e;j)).
    \end{split}
  \end{equation}
  On the other hand,
  \begin{equation}
    \label{eq:compression_p9}
    \weight(\mu;e)=\prod_{j=1}^\infty w\left(\lambda(e;j),\theta(e;j)\right)
  \end{equation}
  and so~(\ref{eq:compression_s1}) follows by taking the quotient
  of~(\ref{eq:compression_p8}) and (\ref{eq:compression_p9}).
\end{proof}
\def\lset#1,#2.{\setbox2=\hbox{$#2$\unskip}F(#1;\ifdim\wd2>0pt
  #2\else k\fi)}
\begin{cor}
  \label{cor:compression}
  Suppose that $W_0$ satisfies the assumptions of
  Theorem~\ref{thm:compression} and let $p\in G$. Assume that for some
  $C_0>0$ one has:
  \begin{equation}
    \label{eq:cor_compression_s1}
    \distance(p,\spt\Gamma)\approx_{C_0}\sigma_k.
  \end{equation}
  Then there is a $C_1=C_1(C_0,\jcut)$ such that:
  \begin{equation}
    \label{eq:cor_compression_s2}
    \lpnorm{\frac{d\xpect[\|\Gamma\|]}{d\mu_p}},,.\approx_{C_1}\sum_{l=1}^k(w_{\gsym}^{-1}\cgwa
    1.)^{l(Q-1)}\sigma_{k-l}.
  \end{equation}
\end{cor}
\begin{proof}
  By assumption~(\ref{eq:cor_compression_s1}) we have that on the
  edges of $\spt\Gamma$:
  \begin{equation}
    \label{eq:ccompression_p1}
    \frac{d\mu_p}{d\mu}\approx_{C(C_0)}(\cgwa 1.\cgwa 2.)^{-k}\prod_{n=k+1}^\infty
    w\left((\lambda(p;n),\theta(p;n))
    \right)^{-1}.
  \end{equation}
  We now obtain the following estimate using that
  $W_0|[\tau_{i+1},\len W_0]$ has a number of edges $\lesssim \sigma_i$:
  \begin{equation}
    \label{eq:ccompression_p2}
    \begin{split}
      \lpnorm{\frac{d\xpect[\|\Gamma\|]}{d\mu_p}},,.&=\Biggl(
      \sum_{\substack{e:\\\indx e.\in[0,\tau_{k-1})}}+\sum_{i=k-1}^{1}\sum_{\substack{e:\\\indx
          e.\in[\tau_{i+1},\tau_{i})}}\\&\mskip \munsplit +\sum_{\substack{e:\\\indx e.\in[\tau_1,\len W]}}
      \Biggr)\Biggl(
      \biggl(\frac{d\xpect[\|\Gamma\|]}{d\mu_p}\mid_e\biggr)^Q\frac{d\mu_p}{d\mu}\mid_e\\
      &\mskip \munsplit \times\weight(\mu;e)
      \Biggr)\\
      &\approx\sum_{\substack{e:\\\indx
          e.\in[0,\tau_{k-1})}}(\cgwa 1.\cgwa 2.)^{-k}\prod_{n=k+1}^\infty
      w\left((\lambda(p;n),\theta(p;n)\right)^{-1}\weight(\mu;e)
      \\
      &\mskip \munsplit + \sum_{i=1}^{k-1}\sum_{\substack{e:\\\indx
          e.\in[\tau_i,\tau_{i+1})}}
      (w_{\gsym}^{-1}\cgwa 1.)^{(k-i)Q}(\cgwa 1.\cgwa 2.)^{-k}\\
      &\mskip \munsplit\mskip \munsplit\times \prod_{n=k+1}^\infty
      w\left((\lambda(p;n),\theta(p;n)\right)^{-1}\weight(\mu;e)
      \\&\mskip \munsplit +\sum_{\substack{e:\\\indx
          e.\in[\tau_k,\len W]}}(w_{\gsym}^{-1}\cgwa 1.)^{kQ}
      (\cgwa 1.\cgwa 2.)^{-k}\\&\mskip \munsplit \mskip \munsplit   \times\prod_{n=k+1}^\infty
      w\left((\lambda(p;n),\theta(p;n)\right)^{-1}\weight(\mu;e)\\
      &\approx\sum_{l=1}^k(w_{\gsym}^{-1}\cgrowth)^{l(Q-1)}\sigma_{k-l}.
    \end{split}
  \end{equation}
\end{proof}
In the following theorem we construct a random curve which moves
``parallel'' to a given walk $W$.
\begin{thm}
  \label{thm:parallel_transport}
  Let $W=\{w_0\,e_1\,w_1\cdots
  e_l\,w_l\}$ be a monotone walk joining $p_0$ to $p_1$ where
  $\ord(p_i)=0$. Let $P_i$ denote the canonical probability measure on
  $\lset p_i,.$. 
  \par To each $p'_0\in\lset p_0,.$ we associate a walk $W_{p'_0}$ as follows. We
  let $w'_0=p'_0$. Then, $e'_{i}$ and (hence) $w'_{i+1}$ are determined by $w'_i$ and
  $e'_{i-1}$ as follows. First $\pi(e'_i)=\pi(e_i)$. If $\ord(w'_i)=0$
  or $w'_i$ is not a gluing or a socket point the previous requirement uniquely
  determines $e'_i$. If $w'_i$ is either a gluing or a socket point of order $>k$ we take
  the edge $e'_i$ satisfying the additional requirement
  $\left(\lambda(e'_i;\ord(w'_{i}),\theta(e'_i;\ord(w'_i))\right)=\left(\lambda(e_i;\ord(w_i)),\theta(e_i;\ord(w_i))\right)$. If $w'_i$ is a
  socket point of order $\le k$ then $e'_i$ is determined by the
  additional requirement that $(\lambda_{e'_i},\theta_{e'_i})=(\lambda_{e'_{i-1}},\theta_{e'_{i-1}})$.
  \par Let $\Gamma$ be the random curve determined by choosing $p'_0$
  according to $P_0$ and letting $\Gamma=W_{p'_0}$. Then the following
  holds:
  \begin{description}
  \item[(C1)] $\cpend \Gamma$ has law $P_1$;
  \item[(C2)] $\spt\Gamma\subset\ball \Gamma_W,C\sigma_k.$;
  \item[(C3)] For $e\in\spt\Gamma$ one has:
    \begin{equation}
      \label{eq:parallel_transport}
      \frac{d\xpect[\|\Gamma\|]}{d\mu}|e\approx_{C_1}=(\cgwa 1.\cgwa
      2.)^{-k}\prod_{j=k}^\infty w\left(
        (\lambda_e(j),\theta_e(j))
      \right)^{-1},
    \end{equation}
    where  $C_1$ depends on {\normalfont \textbf{(P1)}--\textbf{(P3)}} and $\Wgset$.
  \end{description}
\end{thm}
\begin{proof}
  Fix $p'_0\in\lset p_0,.$ and let $e_t$ denote the $t$-th edge of $W$
  and $e'_t$ the $t$-th edge of $W_{p'_0}$. One has
  $\pi(e_t)=\pi(e'_t)$; moreover, the choice of behaviour at gluing
  and socket
  points implies that:
  \begin{equation}
    \label{eq:transport_p1}
    (\lambda(e'_t;j),\theta(e'_t;j))=
    \begin{cases}
      (\lambda(p'_0;j),\theta(p'_0;j))&\text{if $j\le k$}\\
      (\lambda(e_t;j),\theta(e_t;j))&\text{if $j>k$}.
    \end{cases}
  \end{equation}
  Thus, for $e\in\spt\Gamma$ there are a unique $t\in\natural$ and a
  unique $p'_0\in\lset p_0,.$ such that $e$ is the $t$-th edge of
  $W_{p'_0}$. We now prove \textbf{(C1)}. Observe that the end point
  $p'_1$ of $W_{p'_0}$ satisfies:
  \begin{align}
    \label{eq:transport_p2}
    \pi(p'_1)&=\pi(p_1)\\
    (\lambda(p'_1;j),\theta(p'_1;j))&=
    \begin{cases}
      (\lambda(p'_0;j),\theta(p'_0;j))&\text{if $j\le k$}\\
      (\lambda(p_1;j),\theta(p_1;j))&\text{otherwise.}
    \end{cases}
  \end{align}
  Then, using the definition of the map $\tau$ in
  Definition~\ref{defn:label_set}, we get $p'_1=\tau(p'_0)$ and so
  \textbf{(C1)} follows.
  \par Statement \textbf{(C2)} is proven like in
  Theorem~\ref{thm:compression}.
  \par We now show statement \textbf{(C3)}. Let $e\in\spt\Gamma$ and
  let $(t,p'_0)$ be the unique pair such that $e$ is the $t$-th edge
  of $W_{p'_0}$. Then:
  \begin{equation}
    \label{eq:transport_p3}
    \weight\left(
      \xpect[\|\Gamma\|];e
    \right)
    = P(p'_0)=(\cgwa 1.\cgwa 2.)^{-k}\prod_{j=1}^kw\left(
      (\lambda(p'_0;j),\theta(p'_0;j))
    \right),
  \end{equation}
  and the result follows dividing~(\ref{eq:transport_p3}) by $\weight(\mu;e)$.
\end{proof}
\begin{cor}
  \label{cor:transport}
  Let $W$ satisfy the assumptions of
  Theorem~\ref{thm:parallel_transport} and let $p\in G$. Assume that
  for some $C_0>0$ one has:
  \begin{equation}
    \label{eq:transport_s1}
    \distance(p,\spt\Gamma)\approx_{C_0}\sigma_k,
  \end{equation}
  and that $\len W\le C_0\sigma_k$. Then there is a $C_1=C_1(C_0)$
  such that:
  \begin{equation}
    \label{eq:transport_s2}
    \lpnorm{\frac{d\xpect[\|\Gamma\|]}{d\mu_p}},,.\lesssim_{C_1}\sigma_k.
  \end{equation}
\end{cor}
\begin{proof}
  By assumption~(\ref{eq:transport_s1}) we have
  \begin{equation}
    \label{eq:cctransport_p1}
    \frac{d\mu_p}{d\mu}\approx_{C(C_0)}(\cgwa 1.\cgwa 2.)^{-k}\prod_{n=k+1}^\infty\left(
      w((\lambda(p;n),\theta(p;n))
    \right)^{-1}.
  \end{equation}
  on the edges of $\spt\Gamma$. Then for $e\in\spt\Gamma$ one has:
  \begin{equation}
    \label{eq:cctransport_p2}
    \frac{d\xpect[\|\Gamma\|] }
    {d\mu_p }\approx 1.
  \end{equation}
  On the other hand, $\len W\lesssim\sigma_k$ and so:
  \begin{equation}
    \label{eq:cctransport_p3}
    \begin{split}
      \lpnorm{\frac{d\xpect[\|\Gamma\|]}{d\mu_p}},,.&=\sum_{t=1}^{\len
        W}
      \sum_{p'_0\in\lset p_0,.}\sum_{\substack{\text{$e$ is
            the}\\\text{$t$-th edge of}\\\text{$W_{p'_0}$}}}\left(\frac{d\xpect[\|\Gamma\|] }
        {d\mu_p }\right)^{Q}
      \\
      &\mskip \munsplit\times \frac{d\mu_p}{d\mu}\weight(\mu;e)
      \\
      &\approx_{C(C_0)}\sum_{t=1}^{\len
        W}
      \sum_{p'_0\in\lset p_0,.}
      (\cgwa 1.\cgwa 2.)^{-k}\prod_{j=1}^kw\left(
        (\lambda(p'_0;j),\theta(p'_0;j)
      \right)
      \\
      &\lesssim\sigma_k.
    \end{split}
  \end{equation}
\end{proof}
\def\ndec{^{\normalfont\text{(new)}}}
\def\odec{^{\normalfont\text{(old)}}}
\def\lset#1,#2.{\setbox2=\hbox{$#2$\unskip}F(#1;\ifdim\wd2>0pt
  #2\else k-\jcut+1\fi)}
\def\ltset#1,#2.{\setbox2=\hbox{$#2$\unskip}F_\Theta(#1;\ifdim\wd2>0pt
  #2\else k-\jcut+1\fi)}
In the following theorem we assume that the walk is monotone
increasing for concreteness; the same result holds if the walk is
monotone decreasing. The goal is to build a random curve which
``expands'' gaining access to new labels. This is needed to get the
estimate~(\ref{eq:pi_modulus_s2}).
\begin{thm}
  \label{thm:expansion}
  Let $W$ be a monotone increasing walk joining $p_0$ to $p_1$ where
  $\ord(p_i)=0$ and $\len W\in[\sigma_k/2,\sigma_k]$. Assume that all
  edges in $W$ have the same label. Then there is a $C_0$
  which depends only on {\normalfont \textbf{(P1)}--\textbf{(P3)}} such that the
  following holds whenever $\jcut\ge C_0$. Let:
  \begin{equation}
    \label{eq:expansion_st-1}
    \left(
      \lambda(p_0;k-\jcut+1), \theta(p_0;k-\jcut+1)
    \right) = (s_0,t_0),
  \end{equation}
  and choose $(s_1,t_1)\in\Symbset_1\times\Symbset_2\setminus\{(s_0,t_0)\}$.
  \par Choose by Lemma~\ref{lem:mon_lab_inc} a monotone increasing walk
  $W_0\ndec$ from $p_0$ to a socket point $\xi$ of order $k-\jcut+1$,
  and which satisfies $\len W_0\ndec\le\len W$. Let
  $\hat p_1\in\lset p_1,.$ be the point satisfying:
  \begin{equation}
    \label{eq:expansion_s0}
    \left(\lambda(\hat p_1;j),\theta(\hat p_1;j)\right)=
    \begin{cases}
      \left(\lambda(p_1;j),\theta(p_1;j)\right)&\text{for $j\ne k-\jcut+1$}\\
      (s_1,t_1)&\text{for $j=k-\jcut+1$.}
    \end{cases}
  \end{equation}
  Using
  Lemma~\ref{lem:mon_lab_dec} obtain a monotone 
  increasing walk
  $W_{1/2}\ndec$ from $\xi$ to a point $\hat p_{1/2}$ such that
  $\lambda_{\hat p_{1/2}}=\lambda_{\hat p_1}$, $\theta_{\hat
    p_{1/2}}=\theta_{\hat p_1}$ and:
  \begin{equation}
    \label{eq:expansion_s-1}
    \len W_{1/2}\ndec\le\len W-\len W_0\ndec. 
  \end{equation}
  Finally concatenate $W_{1/2}\ndec$ with a monotone increasing walk
  whose edges have constant label $(\lambda_{\hat p_1},\theta_{\hat p_1})$ to obtain a
  walk $W_1\ndec$ joining $\xi$ to $\hat p_1$ and satisfying:
  \begin{equation}
    \label{eq:expansion_s-2}
    \len W_0\ndec + \len W_1\ndec = \len W.
  \end{equation}
  \par Construct a random curve as follows. Choose $p'_0\in\lset
  p_0,k-\jcut.$ using the probability $P_0$. Then with probability:
  \begin{equation}
    \label{eq:expansion_s2}
    (\cgwa 1.\cgwa 2.)^{-1}w((s_0,t_0))\quad(\text{\normalfont event $E\odec$})
  \end{equation}
  let $\Gamma$ be the canonical path $\Gamma_{p'_0\cdot W}$ associated
  to $p'_0\cdot W$. For 
  $(s,t)\ne (s_0,t_0)$ let $\hat
  p'_{1,s,t}$ be the point in $\lset p_1,.$ such that:
  \begin{equation}
    \label{eq:expansion_s-3}
    \left(    \lambda(\hat p'_{1,s,t};j),\theta(\hat p'_{1,s,t};j)
    \right)
    =
    \begin{cases}
      \left(\lambda(p'_0;j),\theta(p'_0;j)\right)&\text{if $j\ne k-\jcut+1$}\\
      (s,t)&\text{if $j=k-\jcut+1$.}
    \end{cases}
  \end{equation}
  Then with 
  probability:
  \begin{equation}
    \label{eq:expansion_s3}
    (\cgwa 1.\cgwa 2.)^{-1}w\left((s,t)\right)\quad(\text{\normalfont event $E_{s,t}\ndec$}),
  \end{equation} let $\Gamma$ be the 
  canonical path associated to the walk:
  \begin{equation}
    \label{eq:expansion_s-4}
    p'_0\cdot W_0\ndec * \left(
      \hat p'_{1,s,t}\cdot(W_1\ndec )^{-1}
    \right)^{-1}.
  \end{equation}
  Then the following hold:
  \begin{description}
  \item[(C1)] $\cpend\Gamma$ has law $P_1$ on $\lset p_1,.$;
  \item[(C2)] $\spt\Gamma\subset\ball\Gamma_W,C\sigma_{k-\jcut+1}.$;
  \item[(C3)] Let
    \begin{equation}
      \label{eq:expansion_s-5}
      E\ndec = \bigcup_{(s,t)\ne(s_0,t_0)}E_{(s,t)}\ndec;
    \end{equation}
    let $\Gamma\odec$ denote $\Gamma$ conditioned on $E\odec$ and
    $\Gamma\ndec$ denote $\Gamma$ conditioned on $E\ndec$. Then for
    each $e\in\spt\Gamma$ there is a unique $\indx e.\in\natural$ such that
    $\pi(e)=\pi(e_{\indx e.})$ where $e_{\indx e.}$ is the $\indx
    e.$-th edge of $W$. If $e\in\spt\Gamma\ndec$ one has:
    \begin{equation}
      \label{eq:expansion_s4}
      \frac{d\xpect[\|\Gamma\ndec\|]}{d\mu}|e\approx_{C_1}
      \cgwa 1.^{-T(e)}\,w_{\gsym}^{-k+T(e)}\cgwa 2.^{-k}\prod_{j=k}^\infty
      w\left((\lambda(e;j),\theta(e;j))\right)^{-1},
    \end{equation}
    where
    \begin{equation}
      \label{eq:expansion_s4bis}
      T(e)=    \begin{cases}      \lg(\len W_0\ndec-\indx e.)
        &\text{if $\max(\pi(e))\le\pi(\xi)$}\\
        \lg(\indx e.-\len W_0\ndec)&\text{otherwise};
      \end{cases}
    \end{equation}
    and if $e\in\spt\Gamma\odec$ then:
    \begin{equation}
      \label{eq:expansion_s5}
      \frac{d\xpect[\|\Gamma\odec\|]}{d\mu}|e\approx_{C_1}(\cgwa
      1.\cgwa 2.)^{-k}\prod_{j=k}^\infty w\left((\lambda(e;j),\theta(e;j))\right)^{-1},
    \end{equation}
    where  $C_1$ depends on $\jcut$, {\normalfont \textbf{(P1)}--\textbf{(P3)}} and $\Wgset$.
  \end{description}
\end{thm}
\begin{rem}
  \label{rem:tatiana_toro_stupid_4}
  Theorem~\ref{thm:expansion} corresponds to the notion of
  ``expanding'' pencils of curves as discussed by Heinonen and 
  Semmes~\cite{semmes_finding_curves, heinonen_analysis}. However,
  here there is a substantial difference with previously known
  examples of PI-spaces, as we need to pass through a socket point in
  order to expand the random curve (or the pencil). This process
  entails some degree of ``compression'' in the expansion, and this
  compression must be controlled as it obstructs the
  Poincar\'e inequality.
  \par Concretely, we want $\Gamma$ to start in $\lset p_0,k-\jcut.$ and
  end in $\lset p_1,.$, where $\jcut$ is an integer parameter chosen
  for convenience, i.e.~to create some ``space'' between the length of $W$
  and the maximum order of entries of $\lambda$ or $\theta$ which
  differ from the corresponding ones in $\lambda_{p_0}$ and $\theta_{p_0}$. While to reach $\lset p_1,
  k-\jcut.$ we can just use a ``parallel lift'' (compare the
  definition of $\Gamma_{p'_0\cdot W}$ using $p'_0\cdot W$), to access
  points $\tilde p_1\in\lset p_1,.$ with $(\lambda_{\tilde
    p_1}(k-\jcut+1),\theta_{\tilde
    p_1}(k-\jcut+1))\ne(\lambda_{p_0}(k-\jcut+1),
  \theta_{p_0}(k-\jcut+1))$ we will
  use the socket point $\xi$.
  \par Specifically, we build a path $W_0\ndec * W_1\ndec$ so that we
  reach from $p_0$ the point $\hat p_1$ whose label is defined
  in~(\ref{eq:expansion_s0}). In this way we can modify the
  $(k-\jcut+1)$-th entry of labels. This construction is then
  generalized to an arbitrary starting point $p'_0\in \lset
  p_0,k-\jcut.$ by using~(\ref{eq:expansion_s-4}).
  \par Heuristically, the event $E\odec$ means that we just follow a
  path ending in $\lset p_1,k-\jcut.$ while the event $E\ndec$ means
  that we pass through $\xi$. Then then technical part of the argument 
  boils down in showing that if the probability of $E\odec$ is chosen
  correctly one gets the estimates~(\ref{eq:expansion_s4})
  and~(\ref{eq:expansion_s5}) which will be needed in verifying the
  Poincar\'e inequality.
\end{rem}
\begin{proof}
  We first explain why the construction of the walks $W_0\ndec$,
  $W_{1/2}\ndec$ and $W_{1}\ndec$ can be carried out. If $C_0$ is
  sufficiently large, one can ensure that whenever $\jcut\ge C_0$, and
  if $C$ is the constant appearing in Lemmas~\ref{lem:mon_lab_inc},
  \ref{lem:mon_lab_dec}, one has:
  \begin{equation}
    \label{eq:expansion_p1}
    2C\sigma_{k-\jcut}\le\len W,
  \end{equation}
  and thus one can construct $W_0\ndec$ and $W_{1/2}\ndec$ satisfying:
  \begin{equation}
    \label{eq:expansion_p2}
    \len W_0\ndec + \len W_{1/2}\ndec \le \len W.
  \end{equation}
  \par We now explain why the concatenation
  in~(\ref{eq:expansion_s-4}) is well-defined. Note that $W_0\ndec$
  and $(W_1\ndec)^{-1}$ satisfy the assumptions of
  Theorem~\ref{thm:compression}; referring to the notation of
  Theorem~\ref{thm:compression}, we have to set $K=k$ where $k$ is now
  given by the integer $k-\jcut+1$
  used in this Theorem; for $W_0\ndec$ the value of $\jcut$ now used in
  Theorem~\ref{thm:compression} is $0$, while for $(W_1\ndec)^{-1}$
  the value of $\jcut$ now used in Theorem~\ref{thm:compression} is $1$. Now, Theorem~\ref{thm:compression}
  ensures that both $p'_0\cdot W_0\ndec$ and $(\hat
  p'_{1,s,t}\cdot W_1\ndec)^{-1}$ end at the point $\xi'\in\ltset \xi,.$
  such that $\theta(p'_0;l)=\theta(\xi';l)$ for $l\ne k-\jcut+1$.
  Therefore, the concatenation
  in~(\ref{eq:expansion_s-4}) is well-defined.
  \par We now turn to the proof of \textbf{(C1)}. Let
  $p'_0=\cpstart\Gamma$; conditional on the event $E\odec$ one has
  that $\cpend\Gamma=p'_1$ where $p'_1$ is the point of $\lset p_1,.$
  satisying $(\lambda_{p'_0},\theta_{p'_0})=(\lambda_{p'_1},\theta_{p'_1})$. The probability of the
  event:
  \begin{equation}
    \label{eq:expansion_p3}
    \left\{
      \cpstart\Gamma=p'_0
    \right\} \cap
    E\odec
  \end{equation}
  is:
  \begin{multline}
    \label{eq:expansion_p4}
    (\cgwa 1.\cgwa 2.)^{-k+\jcut}\prod_{n=1}^{k-\jcut}w\left(
      (\lambda(p'_0;n),\theta(p'_0;n))
    \right) \cdot(\cgwa 1.\cgwa 2.)^{-1} w\left(
      (s_0,t_0)
    \right)\\
    = (\cgwa 1.\cgwa 2.)^{-k+\jcut-1}\prod_{n=1}^{k-\jcut+1}w\left(
      (\lambda(p'_0;n),\theta(p'_0;n))
    \right).
  \end{multline}
  Conditional on the event $E_{s,t}\ndec$ one has $\cpend\Gamma=\hat
  p'_{1,s,t}$, and the probability of the event
  \begin{equation}
    \label{eq:expansion_p5}
    \left\{
      \cpstart\Gamma=p'_0
    \right\} \cap
    E_{s,t}\ndec
  \end{equation} is:
  \begin{multline}
    \label{eq:expansion_p6}
    (\cgwa 1.\cgwa 2.)^{-k+\jcut}\prod_{n=1}^{k-\jcut}w\left(
      (\lambda(p'_0;n),\theta(p'_0;n))
    \right) \cdot(\cgwa 1.\cgwa 2.)^{-1} w\left((s,t)\right)\\
    =(\cgwa 1.\cgwa 2.)^{-k+\jcut-1}\prod_{n=1}^{k-\jcut+1}w\left(
      (\lambda(\hat p'_{1,s,t};n),\theta(\hat p'_{1,s,t};n))
    \right).
  \end{multline}
  We thus conclude that \textbf{(C1)} holds
  \par For \textbf{(C2)} we can apply the same argument as in
  Theorem~\ref{thm:compression}.
  \par We now prove \textbf{(C3)}. The fact that $\indx e.$ is
  well-defined follows from the monotonicity of the walks
  $W$, $W_0\ndec$ and $W_1\ndec$. As all edges of $W$ have the
  same label, for $p'_0\in\lset p_0,k-\jcut.$ one has that
  $p'_0\cdot W = W_{p'_0}$, where $W_{p'_0}$ is defined as in
  Theorem~\ref{thm:parallel_transport}. Therefore, the
  estimate~(\ref{eq:expansion_s5}) on the Radon-Nikodym derivative
  of $\xpect\left[ \|\Gamma\odec\| \right]$ can be obtained
  from~(\ref{eq:parallel_transport}). Let now $t_\xi=\indx e_\xi.$
  where $e_\xi$ is the last edge of $W_0\ndec$. As remarked above,
  the walk $W_0\ndec$
  satisfies the assumptions of Theorem~\ref{thm:compression}. Thus, if
  $e\in\spt\Gamma\ndec$ and $\indx e.\le t_\xi$ we can
  apply~(\ref{eq:compression_s1}) to get~(\ref{eq:expansion_s4})
  with $T(e)=\lg(\len W_0 - \indx e.)$. On the other hand, also the path
  $(W_1\ndec)^{-1}$ satisfies the assumptions of
  Theorem~\ref{thm:compression}. In this case the point
  $\cpend\Gamma\ndec$ avoids the sets of points $p'_1\in\lset
  p_1,.$ such that:
  \begin{equation}
    \label{eq:expansion_p7}
    (\lambda(p'_1;k-\jcut+1),\theta(p'_1;k-\jcut+1))=(s_0,t_0);
  \end{equation}
  in applying Theorem~\ref{thm:compression} this can only
  introduce a multiplicative error lying in
  $[(\cgwa 1.\cgwa 2.)^{-1},\cgwa 1.\cgwa 2.]$ in the
  estimate~(\ref{eq:compression_s1}). Note also that if $\indx
  e.\ge t_\xi$, considering the reverse walk $(W_1\ndec)^{-1}$,
  the integer $\indx e.$ in~(\ref{eq:compression_s1}) must be replaced
  with $\len W-\indx e.$ and thus the proof
  of~(\ref{eq:expansion_s4}) is complete.
\end{proof}
\begin{cor}
  \label{cor:expansion}
  Let $W$ be as in Theorem~\ref{thm:expansion} and let $p\in
  G$. Assume that for some $C_1>0$ one has:
  \begin{equation}
    \label{eq:cor_expansion_s1}
    \distance(p,\spt\Gamma)\approx_{C_1}\sigma_k.
  \end{equation}
  Then there is a $C_2=C_2(C_1,\jcut)$ such that:
  \begin{equation}
    \label{eq:cor_expansion_s2}
    \lpnorm{\frac{d\xpect[\|\Gamma\|]}{d\mu_p}},,.\approx_{C_2}\sum_{l=1}^k(w_{\gsym}^{-1}\cgrowth)^{l(Q-1)}\sigma_{k-l}.
  \end{equation}
\end{cor}
\begin{proof}
  We first apply convexity of the $Q$-th power of the $\lpspace \mu_p,
  .$ norm to get:
  \begin{equation}
    \label{eq:cor_expansion_p1}
    \begin{split}
      \lpnorm{\frac{d\xpect[\|\Gamma\|]}{d\mu_p}},,.&=\lpnorm{P(E\ndec)\frac{d\xpect[\|\Gamma\ndec\|]}{d\mu_p}+P(E\odec)\frac{d\xpect[\|\Gamma\odec\|]}{d\mu_p}},,.
      \\
      &\le
      P(E\ndec)\lpnorm{\frac{d\xpect[\|\Gamma\ndec\|]}{d\mu_p}},,.
      + P(E\odec)\lpnorm{\frac{d\xpect[\|\Gamma\odec\|]}{d\mu_p}},,.;
    \end{split}
  \end{equation}
  let $t_\xi=\indx e_\xi.$ where $e_\xi$ is the last edge of
  $W_0\ndec$.  By assumption~(\ref{eq:cor_expansion_s1}) we can apply
  Corollary~\ref{cor:compression} to $\Gamma\ndec|[0,t_\xi]$ and
  $\Gamma\ndec|[t_\xi,\len W]$. Similarly, by
  assumption~(\ref{eq:cor_expansion_s1}) we can apply
  Corollary~\ref{cor:transport} to $\Gamma\odec$. Thus,
  (\ref{eq:cor_expansion_s2}) follows
  substituting~(\ref{eq:cor_compression_s2}), and (\ref{eq:transport_s2})
  in (\ref{eq:cor_expansion_p1}).
\end{proof}
\subsection{Proof of the Poincar\'e inequality}
\label{subsec:piproof}
In this subsection we join the random curves constructed in
Subsection~\ref{subsec:rand_curves} to prove the Poincar\'e inequality.
\begin{defn}
  \label{defn:neck_range}
  Given $P\ge 1$ we denote by $Q$ the conjugate exponent
  $P/(P-1)$. Let $\neckrange$ denote the range of exponents $P\ge 1$
  such that there is a $C=C(P)$ such that for each $k\in\natural$ one
  has:
  \begin{equation}
    \label{eq:neck_range_1}
    \sum_{l=1}^k(w_{\gsym}^{-1}\cgwa
    1.)^{l(Q-1)}\frac{\sigma_{k-l}}{\sigma_k}\le C.
  \end{equation}
  As $m_k\ge 2$ one has that:
  \begin{equation}
    \label{eq:neck_range_2}
    \left(
      \log_2(w_{\gsym}^{-1}\cgwa 1.) +1 , \infty
    \right)
    \subset\neckrange;
  \end{equation}
  on the other hand, if all $m_k$ are equal to some $m$ one has:
  \begin{equation}
    \label{eq:neck_range_3}
    \neckrange =
    \left(
      \log_m(w_{\gsym}^{-1}\cgwa 1.) +1 , \infty
    \right).
  \end{equation}
\end{defn}
\begin{thm}
  \label{thm:poinc_proof}
  For $P\in\neckrange$ the metric measure space $(G,\mu)$ satisfies a $(1,P)$-Poincar\'e
  inequality, i.e.~$\neckrange\subset\pirange(G,\mu)$.
\end{thm}
\begin{proof}\def\edec{^{\normalfont\text{(exp)}}}
  \def\kdec{^{\normalfont\text{(neck)}}}
  \def\lset#1,#2.{\setbox2=\hbox{$#2$\unskip}F(#1;\ifdim\wd2>0pt
    #2\else k-\jcut\fi)}
  \def\ltset#1,#2.{\setbox2=\hbox{$#2$\unskip}F_\Theta(#1;\ifdim\wd2>0pt
    #2\else k-\jcut\fi)}
  \def\lpnorm#1,#2,#3.{{\setbox0=\hbox{$#2$}\setbox1=\hbox{$#3$}\left\|#1\right\|_{\text{\normalfont L}^{\ifdim\wd0>0pt
          #2 \else Q\fi}(\ifdim\wd1>0pt #3\else \mu_x\fi)}^{\ifdim\wd0>0pt
        #2 \else Q\fi}}}
  We apply Theorem~\ref{thm:pi_modulus}, i.e.~for any pair of points
  $(x,y)$ we show the existence of a random curve $\Gamma$ satisfying:
  \begin{align}
    \label{eq:poinc_proof_p1}
    \spt\Gamma&\subset B(\{x,y\},Cd(x,y)),\\
    \label{eq:poinc_proof_p2}
    \lpnorm{\frac{d\xpect[\|\Gamma\|]}{d(\mu_x+\mu_y)}},,{\mu_x+\mu_y}. &\lesssim_{C_Q}d(x,y),
  \end{align}
  where $C$ does not depend on $x,y$, and $C_Q$ does not depend on
  $x,y$ but depends on $Q$. $\Gamma$ is built by concatenating
  curves obtained by using Theorems~\ref{thm:compression},
  \ref{thm:parallel_transport}, \ref{thm:expansion}. We observe
  that if $\cpend \Gamma_0=\cpstart\Gamma_1$ the random curves
  $\Gamma_0,\Gamma_1$, up to translating their domains, can be
  concatenated to obtain a random curve $\Gamma_0*\Gamma_1$.
  \vskip10pt
  \par\noindent\texttt{Step 1: First part of building ``half'' of a
    random curve joining $x$ to $y$. }
  \par Fix points $x,y$ and assume that $\max\natural(x,y)\le\lg
  d(x,y)$. This assumption will be removed in \texttt{Step 2}. Using
  Theorem~\ref{thm:gw_part1} we can choose a good walk from $x$ to
  $y$ satisfying \textbf{(GWA1)} and \textbf{(GWA2)}. We let $K=\lg
  d(x,y)$. We thus have a uniform constant $C_0$ such that:
  \begin{align}
    \label{eq:poinc_proof_p3}
    C_0\sigma_K&\ge\len W\\
    d(x,w_i)&\ge C_0^{-1}i\quad(\text{$w_i\in W$ is the $i$-th vertex}).
  \end{align}
  For the moment let $C$ be the maximum of the constants occurring
  at points \textbf{(C2)} of Theorems~\ref{thm:compression},
  \ref{thm:parallel_transport}, \ref{thm:expansion}. We can find
  $C_1=C(C_0)$, $J_1=J(C_0)$ such that, if $J\ge J_1$ and $\tilde
  w$ satisfies:
  \begin{equation}
    \label{eq:poinc_proof_p4}
    d(\tilde w,w_i)\le C\sigma_{\lg i -J},
  \end{equation}
  then one has:
  \begin{equation}
    \label{eq:poinc_proof_p5}
    d(\tilde w,x)\ge C_1^{-1}i.
  \end{equation}
  \par We now subdivide $W$ into subwalks $\{W_\alpha\}_{\alpha\in
    I}$ ($I$ is a finite set of integers), the idea
  being that $W$ can be thought of as a concatenation of the
  $\{W_\alpha\}$. More precisely, this can be formalized by using a
  strictly increasing map $\alpha\mapsto m_\alpha$, and letting
  $W_\alpha$ denote the part of $W$ starting at the $m_\alpha$-th
  vertex $w_{m_\alpha}$ and ending at the $m_{\alpha+1}$-th vertex
  $w_{m_{\alpha+1}}$. Note that we obtain an order relation $<$ on $\{W_\alpha\}_{\alpha\in
    I}$ where $W_\alpha<W_{\alpha+1}$.
  \par Using the properties of the good walk constructed in
  Theorem~\ref{thm:gw_part1} we obtain a $J_2$ such that there is a
  decomposition of  $W$ into monotone subwalks $\{W_\alpha\}_{\alpha\in
    I}$ having the following properties:
  \begin{description}
  \item[(Dec1)] For each $k\in\{J_2,\cdots, K\}$ there is a
    $W_\alpha=W\edec_k$ satisfying the assumptions of
    Theorem~\ref{thm:expansion} and:
    \begin{equation}
      \label{eq:poinc_proof_p6}
      \distance(W_\alpha,x)\approx_C\sigma_k;
    \end{equation}
  \item[(Dec2)] For each $k\in\natural(x,y)$ such that $\theta_x(k)\ne\theta_y(k)$, there is a
    $W_\alpha=W\kdec_k$ which can be decomposed into subwalks
    $\tilde W_0$, $\tilde W_1$ which satisfy the following: one has
    $\cpend\tilde W_0=w_{s(k)}=\cpstart\tilde W_1$; moreover, for
    $\jcut\ge J_2$ the walks $\tilde W_0$ and $\tilde W_1^{-1}$
    satisfy the assumptions of Theorem~\ref{thm:compression} where
    $\xi=w_{s(k)}$;
  \item[(Dec3)] For each of the remaining walks $W_\alpha$ there is
    a $k$ such that:
    \begin{align}
      \label{eq:poinc_proof_p7}
      \len W_\alpha&\le C\sigma_k\\
      \label{eq:poinc_proof_p8}
      \distance(W_\alpha,x)&\ge C^{-1}\sigma_k.
    \end{align}       
  \end{description}
  \par $\Gamma$ is constructed by concatenating curves
  $\Gamma_\alpha$ for each $\alpha\in I$. This is done
  inductively, and one starts by letting $\Gamma_1=\Gamma_{W_1}$
  with probability $1$. The next step depends on which of the
  conditions \textbf{(Dec)} is satisfied by $W_{\alpha+1}$:
  \begin{itemize}
  \item Case of \textbf{(Dec1)}. We have $W_{\alpha+1}=W\edec_k$
    and we know that $\cpend\Gamma_\alpha$ is a random point in
    $\lset w_{m_{\alpha+1}},.$ whose law is the canonical
    probability. We obtain $\Gamma_{\alpha+1}$ applying
    Theorem~\ref{thm:expansion}, so that
    $\cpend\Gamma_{\alpha+1}$ is a random point in $\lset
    w_{m_{\alpha+1}},k-\jcut+1.$ whose law is the canonical
    probability. Moreover, by~(\ref{eq:poinc_proof_p6}) we can
    apply Corollary~\ref{cor:expansion} to conclude that:
    \begin{equation}
      \label{eq:poinc_proof_p9}
      \lpnorm{\frac{d\xpect[\|\Gamma_{\alpha+1}\|]}{d\mu_x}},,.\approx_{C_2}\sum_{l=1}^k(w_{\gsym}^{-1}\cgrowth)^{l(Q-1)}\sigma_l,
    \end{equation}
    where $C_2$ is a uniform constant depending on the constants
    $C_0,C_1,C,J_0,J_1,\jcut$. Moreover, by the assumption on $P$
    we have that there is a uniform constant $C_3$ depending on
    $C_2$ and $Q$ such that:
    \begin{equation}
      \label{eq:poinc_proof_p10}
      \lpnorm{\frac{d\xpect[\|\Gamma\|]}{d\mu_x}},,.\lesssim_{C_3}\sigma_k.
    \end{equation}
  \item  Case of \textbf{(Dec2)}. We have $W_{\alpha+1}=W\kdec_k$
    and we know that $\cpend\Gamma_\alpha$ is a random point in
    $\ltset w_{m_{\alpha+1}},.$ whose law is the canonical
    probability. We apply Theorem~\ref{thm:compression} to build
    $\tilde\Gamma_0$ from $\tilde W_0$. We then take the canonical
    probability on $\ltset w_{m_{\alpha+2}},.$ and use
    again Theorem~\ref{thm:compression} to build
    $\tilde\Gamma_1$ from $\tilde W_1^{-1}$. We obtain
    $\Gamma_{\alpha+1}$ by concatenating $\tilde\Gamma_0$ and
    $\tilde\Gamma_1^{-1}$ subject to the following additional
    prescription; suppose that $\cpstart\tilde\Gamma_0=p'_0$;
    then one takes $\cpstart\tilde\Gamma_1=\tau(p'_0)$ where
    $\tau:\lset w_{m_{\alpha+1}},.\to\lset
    w_{m_{\alpha+2}},.$ is the canonical map of
    Definition~\ref{defn:label_set}. Note that:
    \begin{equation}
      \label{eq:poinc_proof_p11}
      \spt\tilde\Gamma_0\cap\spt\tilde\Gamma_1=\{w_{s(k)}\},
    \end{equation}
    as the labels of the edges in $\spt\tilde\Gamma_0$ and
    $\spt\tilde\Gamma_1$ have different $k$-th entries. Moreover,
    as $\xi=w_{s(k)}$ and $d(x,w_{s(k)})\approx_C\sigma_k$, we
    can apply Corollary~\ref{cor:compression} to obtain the
    estimate:
    \begin{equation}
      \label{eq:poinc_proof_p12}
      \lpnorm{\frac{d\xpect[\|\Gamma_{\alpha+1}\|]}{d\mu_x}},,.\approx_{C_2}\sum_{l=1}^k(w_{\gsym}^{-1}\cgrowth)^{l(Q-1)}\sigma_{k-l},
    \end{equation}
    where $C_2$ is a uniform constant depending on the constants
    $C_0,C_1,C,J_0,J_1,\jcut$. Moreover, by the assumption on $P$
    we have that there is a uniform constant $C_3$ depending on
    $C_2$ and $Q$ such that:
    \begin{equation}
      \label{eq:poinc_proof_p13}
      \lpnorm{\frac{d\xpect[\|\Gamma\|]}{d\mu_x}},,.\lesssim_{C_3}\sigma_k.
    \end{equation}
  \item  Case of \textbf{(Dec3)}. We know that $\cpend\Gamma_\alpha$ is a random point in
    $\lset w_{m_{\alpha+1}},.$ and that $\len W_{\alpha+1}\le
    C\sigma_k$. We build $\Gamma_{\alpha+1}$ by applying
    Theorem~\ref{thm:parallel_transport}. In particular, the
    assumptions of Corollary~\ref{cor:transport} are also met an
    so we have:
    \begin{equation}
      \label{eq:poinc_proof_p14}
      \lpnorm{\frac{d\xpect[\|\Gamma_{\alpha+1}\|]}{d\mu_x}},,.\lesssim_{C_2}\sigma_k,
    \end{equation}
    where $C_2$ is a uniform constant depending on the constants
    $C_0,C_1,C,J_0$ and $J_1$.
  \end{itemize}
  \par Note that by the choice of $C_1$, if
  $\spt\Gamma_\alpha\cap\spt\Gamma_\beta\ne\emptyset$, then
  $|\alpha-\beta|\le C_4$, where $C_4$ is a uniform constant. We
  thus obtain that:
  \begin{equation}
    \label{eq:poinc_proof_p15}
    \lpnorm{\frac{d\xpect[\|\Gamma\|]}{d\mu_x}},,.\lesssim_{C_Q}d(x,y)
  \end{equation}
  and that for some uniform $C$:
  \begin{equation}
    \label{eq:poinc_proof_p16}
    \spt\Gamma\subset B(x,Cd(x,y)).
  \end{equation}
  \vskip10pt
  \par\noindent\texttt{Step 2: Modifying Step 1 if
    $\max\natural(x,y)>\lg d(x,y)$.}
  \par In this case $W$ is given by Theorem~\ref{thm:gw_part2}.
  If $\theta_x(\kmax)=\theta_y(\kmax)$ the construction can proceed
  as in \texttt{Step 1} because at $u_{\kmax}$ there is no change of
  the $\theta$-label.
  \par We now discuss the modifications for the case $\theta_x(\kmax)\ne\theta_y(\kmax)$.
  We
  first enlarge $W$ at $w_i=u_{\kmax}$ by inserting $4$ subwalks
  $\{\tilde W_i\}_{i=0}^3$ between $w_i$ and $w_{i+1}$. Let $M=\lg
  d(x,u_{\kmax})$, and let $e$ denote the edge of $W$ before
  $u_{\kmax}$. We take $\tilde W_0$ to be a monotone geodesic walk
  whose edges have all the same label $(\lambda_e,\theta_e)$, with $\len\tilde
  W_0=\sigma_M$ and $d(\tilde W_0,x)\ge C_1^{-1}\sigma_M$. For
  $\tilde W_1$ we take $\tilde W_0^{-1}$. Let now $e$ denote the
  edge of $W$ after $u_{\kmax}$. Then $\tilde W_2$ is a monotone
  geodesic walk whose edges have all the same label $(\lambda_e,\theta_e)$, with $\len\tilde
  W_2=\sigma_M$ and $d(\tilde W_2,x)\ge C_1^{-1}\sigma_M$. For
  $\tilde W_3$ we take $\tilde W_2^{-1}$.
  \par One then proceeds as in \texttt{Step 1}, by subdividing
  $W$. The subdivision must satisfy the additional requirement that
  the $\{\tilde W_i\}_{i=0}^3$ are subwalks of the subdivision, and
  we have only to specify how to construct the corresponding
  $\{\tilde\Gamma_i\}_{i=0}^3$. On $\tilde W_0$ we apply
  Theorem~\ref{thm:parallel_transport} and
  Corollary~\ref{cor:transport} and obtain the estimate:
  \begin{equation}
    \label{eq:poinc_proof_p17}
    \lpnorm{\frac{d\xpect[\|\tilde\Gamma_{0}\|]}{d\mu_x}},,.\lesssim_{C_2}\sigma_M.
  \end{equation}
  Then $\tilde\Gamma_1$ and $\tilde\Gamma_2$ are built by applying
  Theorem~\ref{thm:compression} and Corollary~\ref{cor:compression}
  to $\tilde W_1$ and $\tilde W_2^{-1}$ respectively. Note that
  $\cpstart\tilde\Gamma_2$ is taken to be a random point in
  $\lset\cpstart\tilde W_2^{-1},M-\jcut.$ whose law is the canonical
  probability. We build $\tilde\Gamma_{12}$ by concatenating $\tilde\Gamma_1$ and
  $\tilde\Gamma_2^{-1}$ with the additional prescription that if
  $\cpstart\tilde\Gamma_1=p'_1$ then
  $\cpstart\tilde\Gamma_2=\tau(p'_1)$ where
  \begin{equation}
    \tau:\lset\cpstart\tilde W_1,M-\jcut.\to\lset
    \cpstart\tilde W_2^{-1},M-\jcut.
  \end{equation}
  is the canonical map of
  Definition~\ref{defn:label_set}. We thus obtain the estimate:
  \begin{equation}
    \label{eq:poinc_proof_p18}
    \lpnorm{\frac{d\xpect[\|\tilde\Gamma_{12}\|]}{d\mu_x}},,.\lesssim_{C_3}\sigma_M.
  \end{equation}
  Finally, $\tilde\Gamma_3$ is obtained by applying
  Theorem~\ref{thm:parallel_transport} and
  Corollary~\ref{cor:transport} to $\tilde W_3$. We then have the
  estimate:
  \begin{equation}
    \label{eq:poinc_proof_p19}
    \lpnorm{\frac{d\xpect[\|\tilde\Gamma_3\|]}{d\mu_x}},,.\lesssim_{C_2}\sigma_M.
  \end{equation}
  With these modifications, one
  obtains~(\ref{eq:poinc_proof_p15}),~(\ref{eq:poinc_proof_p16})
  where the constants have possibily worsened compared to
  \texttt{Step 1}.
  \vskip10pt
  \par\noindent\texttt{Step 3: building a random curve
    satisfying~(\ref{eq:poinc_proof_p1}),
    (\ref{eq:poinc_proof_p2})}.
  \par  Fix $x,y\in G$ at distance $>1$. We choose a vertex $z$ of
  order $0$ satisfying:
  \begin{align}
    \label{eq:poinc_proof_p20}
    \left|d(z,x)-\frac{d(x,y)}{2}\right|&\le1\\
    \label{eq:poinc_proof_p20bis}
    \left|d(z,y)-\frac{d(x,y)}{2}\right|&\le1;
  \end{align}\def\jcutx{J_{\text{\normalfont cut},x}}
  \def\jcuty{J_{\text{\normalfont cut},y}}
  we then choose $\jcutx$ and $\jcuty$ larger than $J_2$ of
  \texttt{Step(s) 1, 2} such that:
  \begin{align}
    \label{eq:poinc_proof_p21}
    |\jcutx-J_2|&\le 3\\
    \label{eq:poinc_proof_p22}
    |\jcuty-J_2|&\le 3\\
    \label{eq:poinc_proof_p23}
    \lg d(z,x) -\jcutx &=\lg d(z,y)-\jcuty.
  \end{align}
  We then construct random curves $\Gamma_x$ connecting $x$ to
  $\lset z,\lg d(z,x) -\jcutx.$,
  and $\Gamma_y$ connecting $y$ to $\lset z,\lg d(z,y)-\jcuty.$
  using \texttt{Steps 1,2}. Note that~(\ref{eq:poinc_proof_p23})
  implies that $\cpend\Gamma_x$ and $\cpend\Gamma_y$ have the same
  law. We can thus obtain $\Gamma$ by concatenating $\Gamma_x$ and
  $\Gamma_y^{-1}$. Now~(\ref{eq:poinc_proof_p1}) follows
  from~(\ref{eq:poinc_proof_p16}) and~(\ref{eq:poinc_proof_p20}),
  (\ref{eq:poinc_proof_p20bis}). On the other hand,
  (\ref{eq:poinc_proof_p2}) follows
  from~(\ref{eq:poinc_proof_p15}) and:
  \begin{equation}
    \label{eq:poinc_proof_p24}
    \begin{split}
      \lpnorm{\frac{d\xpect[\|\Gamma\|]}{d(\mu_x+\mu_y)}},Q,{{\mu_x+\mu_y}}.
      &\le\lpnorm{\frac{d\xpect[\|\Gamma_x\|]}{d\mu_x}\,\frac{d\mu_x}{d(\mu_x+\mu_y)}+\frac{d\xpect[\|\Gamma_y\|]}{d\mu_y}
        \,\frac{d\mu_y}{d(\mu_x+\mu_y)}},Q,{{\mu_x+\mu_y}}.\\
      &\le
      2^{Q-1}\left(\lpnorm{\frac{d\xpect[\|\Gamma_x\|]}{d\mu_x}},,\mu_x.+\lpnorm{\frac{d\xpect[\|\Gamma_y\|]}{d\mu_y}},,\mu_y.\right)\\
      &\lesssim_{C_Q}d(x,y).
    \end{split}
  \end{equation}
\end{proof}
\subsection{Lack of the Poincar\'e inequality}
\label{subsec:lack_poinc}
To show that a $(1,P)$ Poincar\'e inequality does not hold if $P$ is
sufficiently small, we produce pairs of points such that the modulus
estimate~(\ref{eq:pi_modulus_s1}) does not hold.
\def\bdec{_{\normalfont\text{bad}}}
\begin{lem}
  \label{lem:bad_box}
  Fix a constant $C_0\ge 1$; then there are constants $M=M(C_0),
  l=l(C_0)$ such that the following holds. Let $(\lambda,\theta)$ be
  labels such that $\lambda(j)=\gsym$ for $j\le k+M$. Let
  $m\in\zahlen$ have order $M+k$ and let $R=3C_0\sigma_k$. In the box
  \begin{equation}
    \label{eq:bad_box_s1}
    B\bdec = \bxset {[m-R,m+R]}, {(\lambda,\theta)}, k+l.
  \end{equation}
  select two points $p_0,p_1$ such that:
  \begin{enumerate}
  \item $\pi(p_0)=m-\sigma_k$ and $\pi(p_1)=m+\sigma_k$;
  \item $\lambda_{p_0}=\lambda_{p_1}$;
  \item $\theta(p_0;j)=\theta(p_1;j)$ if $j\ne k+M$ and $\theta(p_0;k+M)\ne\theta(p_1;k+M)$.
  \end{enumerate}
  Then there is a constant $C_1(C_0,P)$ such that:
  \begin{equation}
    \label{eq:bad_box_s2}
    d(p_0,p_1)^{P-1}\modulus , {p_0,p_1},\mu_{p_0,p_1}^{(C_0)}.\le
    \frac{C_1}{(k-1)^P}
    \sum_{i=1}^{k-1}
    \left(
      \frac{\sigma_k}{\sigma_i}
    \right)^{P-1} (w_{\gsym}\cgwa 1.^{-1})^{k-1-i}.
  \end{equation}
\end{lem}
\begin{proof}
  Let $\xi\in B\bdec$ denote the socket point of label
  $(\lambda,\theta)$ such that $\pi(\xi)=m$. Let $\gamma$ be a
  continuous curve joining $p_0$ to $p_1$. Note that by possibly
  enlarging $C_0$ we have $d(p_0,p_1)\approx_{C_0}\sigma_k$ and so for
  $l(C_0)$ sufficienlty large, by Lemma~\ref{lem:ball_box} we have:
  \begin{equation}
    \label{eq:bad_box_p1}
    B(\{p_0,p_1\},C_0d(p_0,p_1))\subset B\bdec.
  \end{equation}
  If $M(C_0)$ is sufficiently large, the only integer of order $k+M$
  contained in $\pi(B\bdec)$ is $m$. To estimate $\modulus ,
  p_0,p_1,\mu_{p_0,p_1}^{(C_0)}.$ we need to produce an appropriate
  Borel function $g$. For the moment we let $g=\infty$ on $B\bdec^c$ and then the case of
  interest becomes when $\gamma$ stays in $B\bdec$; in particular,
  $\gamma$ must pass through a socket point $\xi'\in\ltset\xi,k+l.$.
  \par Let $s\in\dom\gamma$ be the first time when
  $\gamma(s)\in\ltset\xi,k+l.$ and let $\gamma_1=\gamma|[0,s]$. Note
  that:
  \begin{equation}
    \label{eq:bad_box_p2}
    [m-\sigma_{k-1},m]\subset\pi\circ\gamma([0,s]);
  \end{equation}
  for $i<k$ let $t_i=m-\sigma_i$ and let $\varrho_i$ be the last time
  such that $\pi\circ\gamma_1(\varrho_i)=t_i$.
  Let $E(i)$ denote the set of edges $e\in B\bdec$ such that
  $\pi(e)\subset[t_i,t_{i-1}]$. As there are no integers of order $i$
  in $[t_i,m]$ we conclude that the curve
  $\gamma_1|[\varrho_i,\varrho_{i-1}]$ passes through edges
  $\{e_1,\cdots,e_l\}\subset E(i)$ such that:
  \begin{description}
  \item[$(E(i),$1$)$] $e_a$ and $e_{a+1}$ are adjacent,
    $t_i\in\pi(e_1)$ and $t_{i-1}\in\pi(e_l)$;
  \item[$(E(i),$2$)$] $l\ge \sigma_i-\sigma_{i-1}$;
  \item[$(E(i),$3$)$] $\lambda(e_a;j)=\gsym$ for $j\ge i$;
  \item[$(E(i),$4$)$] $\theta(e_a;j)=\theta(p_0;j)$ for $j>k+l$.
  \end{description}
  We now complete the definition of $g$ by defining $g|B\bdec$ as follows: if $e\in E(i)$ for
  some $i$ and \textbf{$(E(i),$3$)$} and \textbf{$(E(i),$4$)$} hold,
  we let $g=(k-1)^{-1}(\sigma_i-\sigma_{i-1})^{-1}$; otherwise, we let
  $g=0$. We now obtain the following lower bound:
  \begin{equation}
    \label{eq:bad_box_p3}
    \begin{split}
      \int g\,d\hmeas
      1._\gamma&\ge\sum_{i=1}^{k-1}\int\chi_{E(i)}g\,d\hmeas 1._\gamma\\
      &\ge\sum_{i=1}^{k-1}\int_{\varrho_i}^{\varrho_{i-1}}\chi_{E(i)}(\gamma(\tau))g(\gamma(\tau))\,d\tau\\
      &\ge\sum_{i=1}^{k-1}\frac{\sigma_i-\sigma_{i-1}}{(k-1)(\sigma_i-\sigma_{i-1})}=1,
    \end{split}
  \end{equation}
  where we let $\sigma_0=0$.
  \par Note now that $\gamma_1|[\varrho_{k-1},s]$ is at distance
  $\approx_{C_2}\sigma_k$ from $p_0,p_1$, where $C_2$ is a uniform
  constant. Therefore, we have:
  \begin{equation}
    \label{eq:bad_box_p4}
    \frac{d\mu_{p_0,p_1}^{(C_0)}}{d\mu}\mid (\gamma_1[\varrho_{k-1},s])
    \approx_{C(C_2)}(\cgwa 1.\cgwa 2.)^{-k-1}\prod_{j>k-1}
    w\left(
      \lambda(p_0;j),\theta(p_0;j)
    \right)^{-1};
  \end{equation}
  note that $g\ne0$ in $B\bdec$ only on $\bigcup_{i=1}^{k-1}E(i)$, and
  let $\tilde E(i)$ denote the set of edges of $E(i)$ satisfying
  \textbf{$(E(i),$3$)$} and \textbf{$(E(i),$4$)$}; as $g$ vanishes on
  $E(i)\setminus \tilde E(i)$, we have for some $C_1(C_0,C_2,P,M,l)$:
  \begin{equation}
    \label{eq:bad_box_p5}
    \begin{split}
      \int
      g^P\,d\mu_{p_0,p_1}^{(C_0)}&\lesssim_{C_1}\sum_{i=1}^{k-1}\frac{1}{(k-1)^P\sigma_i^P}
      (\cgwa 1.\cgwa 2.)^{-k+1}\\
      &\mskip\munsplit \times\sum_{e\in\tilde E(i)}\prod_{j\le
        k-1}w(\lambda_e(j),\theta_e(j))\\
      &\lesssim_{C_1}\sum_{i=1}^{k-1}\frac{\sigma_i}{(k-1)^P\sigma_i^P}(\cgwa
      1.\cgwa 2.)^{-k+1}\\
      &\mskip\munsplit\times w_{\gsym}^{k-1-i}\cgwa 1.^i\cgwa
      2.^{k+l}\\
      &\lesssim_{C_1}\sum_{i=1}^{k-1}\frac{\sigma_i}{(k-1)^P\sigma_i^P}(w_{\gsym}\cgwa 1.^{-1})^{k-1-i},
    \end{split}
  \end{equation}
  from which~(\ref{eq:bad_box_s2}) follows.
\end{proof}
\begin{thm}
  \label{thm:lack_poinc}
  If $P\le 1+\log_N(w_{\gsym}^{-1}\cgwa 1.)$ then $P\not\in\pirange
  (G,\mu)$. Moreover, if all the $m_k$ are equal to some $m$,
  $\neckrange=\pirange(G,\mu)$.
\end{thm}
\begin{proof}
  We show that for any value of $C$, (1) in
  Theorem~\ref{thm:pi_modulus} fails. For any $k\ge 1$ we can find a
  \emph{bad box} $B\bdec$ satisfying the assumptions of
  Lemma~\ref{lem:bad_box}. Hence we find sequences of pairs of points
  $(p_0^{(k)},p_1^{(k)})\in G^2$ such that:
  \begin{equation}
    \label{eq:lack_poinc_p1}
    d\left(
      p_0^{(k)},p_1^{(k)}
    \right)^{P-1}\modulus
    ,{p_0^{(k)},p_1^{(k)}},.\le\frac{C}{(k-1)^P}\sum_{i=1}^{k-1}
    \left(
      \frac{\sigma_k}{\sigma_i}
    \right)^{P-1}(w_{\gsym}\cgwa 1.^{-1})^{k-1-i}.
  \end{equation}
  As $P\le 1+\log_N(w_{\gsym}^{-1}\cgwa 1.)$ and as
  $\sigma_k/\sigma_i\le N^{k-i}$, the
  rhs.~of~(\ref{eq:lack_poinc_p1}) goes to $0$ as
  $k\nearrow\infty$.
  \par If all the $m_k$ are equal to some $m$, then the
  rhs.~of~(\ref{eq:lack_poinc_p1}) goes to $0$ exactly when $P\le
  1+\log_m(w_{\gsym}^{-1}\cgwa 1.)$. 
\end{proof}
\begin{rem}
  \label{rem:arbitrary_range}
  Note that as $w_{\gsym}^{-1}\cgwa 1.\nearrow\infty$ one has $\min
  \pirange(G,\mu)\to \infty$, i.e.~the range of exponents for which a
  Poincar\'e inequality holds gets narrower and narrower. On the other
  hand, as $w_{\gsym}^{-1}\cgwa 1.\searrow1$, $\min\pirange(G,\mu)\to
  1$ and thus the range of exponents for which the Poincar\'e
  inequality holds can be arbitrarily prescribed. However, as either
  $w_{\gsym}^{-1}\cgwa 1.$ goes to $1$ or $\infty$, the doubling
  constant of $\mu_G$ blows up.
\end{rem}
\section{Stability under blow-up}
\label{sec:blowup}
In this section we show how to use $G$ to construct a metric measure
space $X$ which satisfies the conclusion of Theorem~\ref{thm:main_result}. From
now on the measure on $G$ that we constructed in Section~\ref{sec:exa}
will be denoted by $\mu_G$. In this section we often deal with balls
of different spaces, and so at times we add a subscript to them to
distinguish the space to which they belong. Given a metric
space $X$, we use $\lambda X$ to denote $X$ with the metric rescaled
by the factor $\lambda>0$.
\subsection{Asymptotic cones}
\label{subsec:acones}
In this subsection we define asymptotic cones and construct the
example $X$.
\begin{defn}
  \label{defn:asy_cone}
  An \textbf{asymptotic cone} of a metric measure space $(X,\mu)$ is a
  measured pointed Gromov-Hausdorff limit of a sequence of
  rescalings:
  \begin{equation}
    \label{eq:asy_cone_1}
    \left(
      \lambda^{-1}_nX,\frac{\mu }{\mu\left(\sball X,p_n,\lambda_n.\right) },p_n
    \right) 
  \end{equation}
  where $\lim_{n\to\infty}\lambda_n=\infty$. Note that $\sball
  X,p_n,\lambda_n.$ denotes a ball of radius $\lambda_n$
  in $X$, that is a ball of radius $1$ in $\lambda^{-1}_nX$. The set of
  asymptotic cones of $(X,\mu)$ will be denoted by $\asycone
  X,\mu.$. Note that it would be more appropriate to say that
  $\asycone X,\mu.$ is a set of equivalence classes of metric
  spaces under measure-preserving isometries, but we will avoid
  such subtleties in the following discussion.
\end{defn}
\begin{defn}
  \label{defn:weak_tan}
  A \textbf{weak tangent} $(Y,\nu,q)$ of a metric measure space $(X,\mu_X)$ is a
  measured pointed Gromov-Hausdorff limit of a sequence of rescalings:
  \begin{equation}
    \label{eq:weak_tan_1}
    \left(
      \lambda_nX,\frac{\mu_X }{\mu_X\left(\sball X,p_n,\lambda_n^{-1}.\right) },p_n
    \right) 
  \end{equation}
  where $\lim_{n\to\infty}\lambda_n=\infty$. The set of weak tangents of $(X,\mu_X)$
  will be denoted by $\wtan ,.$.
\end{defn}
In the case of $(G,\mu_G)$ the fact that asymptotic cones exist and
that the corresponding measures are doubling with uniformly bounded
doubling constants follows from a standard compactness argument.
\begin{lem}
  \label{lem:cone_set_closure}
  The set of asymptotic cones $\asycone ,.$ is closed under the
  operation of taking weak tangents, i.e.~whenever
  $(X,\mu_X,p)\in\asycone ,.$ one has $\wtan ,.\subset\asycone ,.$.
\end{lem}
\begin{proof}
  On the metric level, the proof is straighforward using that one can
  approximate a weak tangent $(Y,\mu_Y,q)\in\wtan ,.$ by rescaling an
  approximating sequence for $(X,\mu_X,p)$. There is, however, an
  issue with normalization of balls which is addressed in the
  following lemma.
\end{proof}
\begin{lem}
  \label{lem:continuity_balls}
  Let $(X,\mu,p)\in\asycone ,.$ and consider a sequence of rescalings:
  \begin{equation}
    \label{eq:continuity_balls_s1}
    \left(
      \lambda_n^{-1}G ,
      \underbrace {
        \frac{ \mu_G }
        {\mu_G\left( \sball X,p_n,\lambda_n. \right)}
      }_{\nu_n},
      p_n
    \right) \to (X,\mu,p). 
  \end{equation}
  Then for each $t\ge 0$ one has:
  \begin{equation}
    \label{eq:continuity_balls_s2}
    \lim_{n\to\infty}\nu_n\left(
      \sball G,p_n,\lambda_nt.
    \right) = \mu \left(
      \sball X,p,t.
    \right).
  \end{equation}
\end{lem}
\begin{proof}
  Using that $n\mapsto\nu_n(G)$ is lower semicontinuous if $G$ is open and
  upper semicontinuous if $G$ is compact, it suffices to show that
  one has, uniformly in $p_n,\lambda_n$:
  \begin{equation}
    \label{eq:continuity_balls_p1}
    \frac{
      \mu_G \left(
        \ball p_n,\lambda_n t. \setminus \ball p_n,\lambda_n(t-\varepsilon).
        \right)
    }
    {
      \mu_G \left(
        \ball p_n,\lambda_nt.
        \right)
      } \le O(\varepsilon^{1/2}).
    \end{equation}
    For $s\in(0,1)$ let $L(s)$ denote the set of labels
    $(\lambda,\theta)$ of edges intersecting $\partial \ball
    p_n,\lambda_n(1-s)t.$. Note that $s_1<s_2$ implies $L(s_2)\supset
    L(s_1)$. However, as:
    \begin{equation}
      \label{eq:continuity_balls_p2}
      \frac{
        \lambda_n(1-s)t
      }
      {
        \lambda_n(1-\varepsilon^{1/2})t
        } \le \frac{3}{2}
      \end{equation}
      for $\varepsilon$ sufficiently small and $s\ge
      \varepsilon$, any label $(\lambda,\theta)\in L(s)\setminus
      L(\varepsilon^{1/2})$ can differ from a label in
      $L(\varepsilon^{1/2})$ only at the $j$-th entry, where either:
      \begin{equation}
        \label{eq:continuity_balls_ins_p1}
        j\in\left\{
          \lg
      (2\lambda_n(1-\varepsilon^{1/2})t), \lg
      (2\lambda_n(1-\varepsilon^{1/2})t)+1
        \right\},
      \end{equation}
      or $j=j_0$, where $j_0$ is some fixed integer
      $>\lg
      (2\lambda_n(1-\varepsilon^{1/2})t)+1$ (this can occur if the
      ball $\ball p_n,\lambda_nt.$ contains a socket point of order
      greater than $\lg (2\lambda_nt)$). We thus obtain:
      \begin{equation}
        \label{eq:continuity_balls_p3}
        \frac{
      \mu_G \left(
        \ball p_n,\lambda_n t. \setminus \ball p_n,\lambda_n(t-\varepsilon).
        \right)
    }
    {
      \mu_G \left(
        \ball p_n,\lambda_nt.
        \right)
      } \le (\cgwa 1.\cgwa 2.)^3 \frac{\lambda_n\varepsilon t}{\lambda_n\varepsilon^{1/2}t},
    \end{equation}
    from which~(\ref{eq:continuity_balls_p1}) follows.
\end{proof}
\par We will use a discretization procedure of Gill and Lopez~\cite{gill_lop_disc} that allows to compare PI
spaces and graphs. We rephrase their result in a
slightly more general context, where there is more freedom in the
choice of the approximating graph; the proof is omitted being a
straightforward generalization of their argument.
\begin{thm}
  \label{thm:gill_disc}
  Let $H$ be a connected graph whose metric is a constant multiple of
  the length metric. For $\varepsilon>0$ and $C_0>0$ consider a subset
  $V$ of vertices of $H$ which is an $\varepsilon$-separated net and
  $C_0\varepsilon$-dense. Assume that for some $C_1>0$ there is a
  $C_1$-biLipschitz embedding $F:V\to X$ such that $F(V)$ is
  $C_1\varepsilon$-dense in $X$. Let $\mu_X$ be a doubling measure on
  $X$ with constant $C_2$. Let $\mu_H$ be a doubling measure on $H$
  which restricts to a multiple of arclength on each edge and such
  that one has, for some $C_3>0$:
  \begin{equation}
    \label{eq:gill_disc_s1}
    \mu_H \left(
      \sball H,v,r.
    \right)\approx_{C_3}
    \mu_X \left(
      \sball X,F(v),r.
    \right) \quad(\forall(v,r)\in V\times[\varepsilon,\infty)).
  \end{equation}
  Then $\pirange (X,\mu_X)\subset \pirange (H,\mu_H)$; moreover, if
  $C_X(P)$ denotes the constant of the $(1,P)$-Poincar\'e inequality
  in $(X,\mu_X)$, then the corresponding constant $C_H(P)$ in $(H,\mu_H)$
  satisfies:
  \begin{equation}
    \label{eq:gill_disc_s2}
    C_H(P)\le C(C_0,C_1,C_2,C_3,C_X(P),\varepsilon).
  \end{equation}
\end{thm}
Since we work with pointed measured Gromov-Hausdorff convergence we
need a \emph{local} version of Theorem~\ref{thm:gill_disc}.
\begin{cor}
  \label{cor:gill_disc}
  In Theorem~\ref{thm:gill_disc}, assume that $V$ is not
  $C_0\varepsilon$-dense in the whole of $H$, but that $V$ now lies
  in a ball $\scball H,h,R.$ with $R>0$ in which it is
  $C_0\varepsilon$-dense. Assume also that $F(V)$ contains a
  $C_1\varepsilon$-dense set in a ball $\sball
  X,x,C_1^{-1}R.$. Furthermore, assume that $X$ is geodesic. Then the
  conclusion of Theorem~\ref{thm:gill_disc} holds replacing
  $(H,\mu_H)$ with:
  \begin{equation}
    \label{eq:cor_gill_disc_s1}
    \left(
      \scball H,h,\tilde{C}^{-1}R. ,
      \mu_H\on\scball H,h,\tilde{C}^{-1}R.
    \right),
  \end{equation}
  where $\tilde{C}$ depends only on $C_0,C_1,C_2$, and $\varepsilon$.
\end{cor}
\begin{proof}
  One can reduce this local case to the \emph{global} one,
  Theorem~\ref{thm:gill_disc}, by recalling that if $(X,\mu)$ is
  geodesic and admits a $(1,P)$-Poincar\'e inequality with exponent
  $C(P)$, there is a $C_1(C(P))$ such that each for each $R>0$ the metric
  measure space $(\clball x,R.,\mu\on\clball x,R.)$ admits a
  $(1,P)$-Poincar\'e inequality with constant $C_1$
  (see \cite{haj_kosk_sobolev_met}).
\end{proof}
\begin{thm}
  \label{thm:bwup_stability}
  Let $(X,\mu_X,p)\in\asycone ,.$; then:
  \begin{equation}
    \label{eq:bwup_stability_s1}
    \pirange(X,\mu_X)=\pirange(G,\mu_G).
  \end{equation}
\end{thm}
\def\ccut{C_{\text{\normalfont cut}}}
\def\cnpi{C_{\text{\normalfont PI}}}
\begin{proof}
  \noindent\texttt{Step 1: $\pirange(X,\mu_X)\subset
    \pirange(G,\mu_G)$.}
  \par Let
  \begin{equation}
    \label{eq:bwup_stability_p1}
    \left(
      \lambda^{-1}_nG,\underbrace{\frac{\mu_G }{\mu_G\left(\sball G,p_n,\lambda_n.\right) }}_{\nu_n},p_n
    \right) \to (X,\mu_X,p)
  \end{equation}
  and assume that $P\in\pirange(X,\mu_X)$, $C(P)$ being the
  corresponding constant. Choose $N(n)$ such that:
  \begin{equation}
    \label{eq:bwup_stability_p2}
    1\le \frac{
      \lambda_n
    }
    {
      \sigma_{N(n)}
    }
    \le N
  \end{equation} and pass to a subsequence such that
  $\lim_{n\to\infty}\frac{\lambda_n}{\sigma_{N(n)}}$
  exists. Therefore, up to rescaling the metric on $X$ by a factor in
  $[1/N,1]$ we can assume that:
  \begin{equation}
    \label{eq:bwup_stability_p3}
    (\underbrace{\sigma^{-1}_{N(n)}G}_{G_n},\nu_n,p_n)\to(X,\mu_X,p);
  \end{equation}
  note also that $(X,\mu_X)$ is geodesic being a limit of geodesic
  metric spaces. Fix $\varepsilon, R>0$; for $n\ge D_0(R,\varepsilon)$
  we can assume that the Gromov-Hausdorff distance between $\sball
  G_n,p_n,R.$ and $\sball X,p,R.$ is at most
  $\frac{\varepsilon}{3}$. Now the vertices of order $\ge l$ in $G$
  form a maximal $\sigma_l$-net which becomes a maximal
  $\sigma_l\sigma^{-1}_{N(n)}$-net in $G_n$; for each $n$ we choose
  $N_\varepsilon(n)\le N(n)$ such that:
  \begin{equation}
    \label{eq:bwup_stability_p4}
    \varepsilon\le \sigma_{N_\varepsilon(n)}\sigma^{-1}_{N(n)}\le N\varepsilon.
  \end{equation}
  Lett $V(n;\varepsilon)$ be the set of vertices of $G_n$ whose order
  in $G$ is at least $N_\varepsilon(n)$ and which are contained in
  $\sball G_n,p_n,R.$. Then $V(n;\varepsilon)$ is an
  $\varepsilon$-separated net in $\sball G_n,p_n,R.$ and is also
  $N\varepsilon$-dense there. Thus the cardinality of
  $V(n;\varepsilon)$ is uniformly bounded in $n$ and
  $V(n;\varepsilon)\to W$ in the Hausdorff sense where $W$ is a
  $\frac{2}{3}\varepsilon$-separated net in $\sball X,p,R.$ in which
  it is also $\frac{3}{2}N\varepsilon$-dense. Therefore for $n\ge
  D_0(R,\varepsilon)$ we find an $L$-biLipschitz map:
  \begin{equation}
    \label{eq:bwup_stability_p5}
    F_n:V(n;\varepsilon)\to W,
  \end{equation}
  where $L$ does not depend on $\varepsilon$ or $n$. Now, as the
  cardinalities of $V(n;\varepsilon)$ and $W$ are uniformly bounded,
  for $n\ge D_1(R,\varepsilon)$ we can assume that the sets
  $V(n;\varepsilon)$ and $W$ have the same cardinality and write
  $V(n;\varepsilon)=\{v_\alpha^{(n)}\}_{\alpha\in A}$ and
  $W=\{w_\alpha\}_{\alpha\in A}$ so that
  $F_n(v_\alpha^{(n)})=w_\alpha$ for each $\alpha\in A$. We now
  use a variation on the argument of Lemma~\ref{lem:continuity_balls}
  (where we take balls not centred on the basepoints) to conclude that for each
  $r\in[\varepsilon,R]$ one has:
  \begin{equation}
    \label{eq:bwup_stability_p6}
    \nu_n\left(
      \sball G_n,v_\alpha^{(n)},R.
    \right) \to
    \mu_X\left(
      \sball X,w_\alpha,R.
    \right);
  \end{equation}
  so for $n\ge D_2(R,\varepsilon)$ we can assume that:
  \begin{equation}
    \label{eq:bwup_stability_p7}
    \nu_n \left(
      \sball G_n,v_\alpha^{(n)},R.
    \right) \approx_{1+\varepsilon}
    \mu_X \left(
      \sball X,w_\alpha,R.
    \right).
  \end{equation}
  We now apply Corollary~\ref{cor:gill_disc} and find
  $\ccut=\ccut(\varepsilon)$ such that
  \begin{equation}
    \label{eq:bwup_stability_p8}
    \left(
      \scball G_n,p_n,R/\ccut.,\nu_n\on\scball G_n,p_n,R/\ccut.
    \right)    
  \end{equation}
  admits a $(1,P)$-Poincar\'e inequality with constant
  $\cnpi=C(C(P),\varepsilon)$. By rescaling back we conclude that:
  \begin{equation}
    \label{eq:bwup_stability_p9}
    \left(
      \scball G,p_n,\sigma_{N(n)}R/\ccut. ,
      \mu_G\on\scball G,p_n,\sigma_{N(n)}R/\ccut.
    \right)
  \end{equation}
  admits a $(1,P)$-Poincar\'e inequality with constant $\cnpi$. Fix a
  basepoint $q\in G$. For each $s>0$ we can find $n\ge D_3(s)$ such
  that $\scball G,p_n,\sigma_{N(n)}R/\ccut.$ contains an isometric
  copy $B_s$ of $\scball G,q,s.$ and such that the measures $\mu_G\on
  B_s$ and $\mu_G\on \scball G,q,s.$ agree up to a multiple. Thus
  \begin{equation}
    \label{eq:bwup_stability_p10}
    \left(\scball G,q,s.,\mu_G\on\ball G,q,s.\right)
  \end{equation}
  admits a $(1,P)$-Poincar\'e inequality with constant $\cnpi$; as
  $\cnpi$ does not depend on $s$ we conclude by letting $s\to\infty$.
  \noindent\par\texttt{Step 2: $\pirange(G,\mu_G)\subset
    \pirange(X,\mu_X)$.}
  \noindent\par This follows from the stability of the Poincar\'e
  inequality under measured pointed Gromov-Hausdorff convergence, see
  \cite{keith-modulus}. 
\end{proof}
\subsection{Putting all together}
\label{subsec:together}
In this subsection we complete the proof of Theorem~\ref{thm:main_result}.
\begin{proof}[Proof of Theorem~\ref{thm:main_result}]
  The existence of the measures $\{\mu_{P_c}\}_{P_c}$ follows
  combining Theorems~\ref{thm:bwup_stability}, \ref{thm:poinc_proof},
  \ref{thm:lack_poinc} and Remark~\ref{rem:arbitrary_range}.
  \par The projection map $\pi:G\to\real$ passes to the limit giving a
  $1$-Lipschitz map $\pi:X\to\real$. The geodesic lines of the form
  $\real\times\{\lambda\}\times\{\theta\}$ pass to the limit and give
  a Fubini-like representation of the measure $\mu_{P_c}$. To this
  Fubini representation one can associate a Weaver derivation $D$, i.e.~a
  horizontal vector field as in~\cite{deralb}.
  \par The verification that $(X,\pi)$ is a chart is standard and can
  be carried out in two ways. The first way uses a Sobolev-space
  argument like Sec.~9 in \cite{cheeger_inverse_poinc}. The second uses $D$ and
  the Stone-Weierstrass Theorem for Lipschitz Algebras as in
  Example\cite[Example~5E]{weaver00}.
  \par The claim about the Assouad-Nagata dimension follows because the graph
  $G$ has Assouad-Nagata dimension $1$ and the Assouad-Nagata dimension is stable in
  passing to asymptotic cones.
\end{proof}
\bibliographystyle{alpha}
\bibliography{neck_poinc_biblio}
\end{document}